\documentclass{amsart}
\usepackage{amsfonts}
\usepackage{amsmath}
\usepackage{amssymb}
\usepackage{mathtools}
\usepackage{graphicx, caption, subcaption}
\usepackage{ulem}
\usepackage{mathrsfs}
\usepackage{tikz-qtree}
\usepackage{tikz-cd}

\usepackage{xcolor}

\address[tlazarus@math.arizona.edu]{Tynan Lazarus, Department of Mathematics, 
 University of Arizona, 617 N. Santa Rita Ave, Tucson, AZ, 85721, USA}
 \address[ealvarado@math.ucdavis.edu]{Enrique Alvarado, Department of Mathematics, University of California at Davis,1 Shields Ave, Davis, CA, USA}
\address[qlxia@math.ucdavis.edu]{Qinglan Xia, Department of Mathematics, University of California at Davis,1 Shields Ave, Davis, CA, USA}

\numberwithin{equation}{section}

\newtheorem{theorem}{Theorem}[section]
\newtheorem{corollary}[theorem]{Corollary}
\newtheorem{lemma}[theorem]{Lemma}
\newtheorem{proposition}[theorem]{Proposition}

\newtheorem{definition}{Definition}[section]
\newtheorem{example}{Example}[section]
\newtheorem{remark}{Remark}[section]

\newcommand{\R}{{\mathbb{R}}}
\newcommand{\N}{{\mathbb{N}}}

\renewcommand{\k}{\mathbf{k}}
\newcommand{\Z}{{\mathbb{Z}}}
\newcommand{\C}{{\mathcal C}}
\newcommand{\F}{{\mathcal F}}
\newcommand{\diam}{\textsl{diam}}
\newcommand{\im}{\mathcal{G}_m}

\newcommand{\overbar}[1]{\mkern 1.5mu\overline{\mkern-1.5mu#1\mkern-1.5mu}\mkern 1.5mu}

\newcommand{\commentedtext}[1]{}
\newcommand{\comments}[1]{}
\newcommand{\blue}[1]{\color{blue} #1 \color{black}}

\newcommand{\mylee}{\,\hbox{\vrule height 0.4pt depth 0pt width 6pt
		\vrule height 7pt depth 0pt}\,}

\makeatletter
\newcommand{\proofstep}[1]{%
  \par
  \addvspace{\medskipamount}
  \textit{#1\@addpunct{.}}\enspace\ignorespaces
}
\makeatother

\title[$\mathcal{F}$-limit sets]{General fractals represented by $\mathcal{F}$-limit sets of compression maps}

\author[T. Lazarus, E. Alvarado, Q. Xia]{Tynan Lazarus, Enrique G Alvarado, Qinglan Xia}

\begin{document}
	
\begin{abstract}
In this article, we provide a simple and systematic way to represent general (inhomogeneous) fractals that may look different at different scales and places. 
By using set-valued compression maps, we express these general fractals as $\mathcal{F}$-limit sets, which are represented as sequences of points in a fixed parameterization space $M$. 
By choosing different types of sequences in $M$, we get
various types of fractals: from self-simlilar to non self-similar, and from deterministic to random.  
The computational complexity of producing a general fractal is independent of the sequence in $M$, and as a result, is the same as that of an iterated function system obtained from a constant sequence.   
In the metric space setting, we also estimate the Hausdorff dimension of limit sets for collections of sets that do not necessarily satisfy the Moran structure conditions.
In particular, we introduce the concept ``uniform covering condition"  for the study of the lower bound of the Hausdorff dimension of the limit set, and provide sufficient conditions for this condition.
Specific examples (Cantor-like sets, Sierpi\'nski-like Triangles, etc.) with the calculations of their corresponding Hausdorff dimensions are also studied. 
\end{abstract}

	
	
\maketitle
\section{Introduction}
\label{intro}
A popular mathematical way to produce a fractal is with the limit set of a collection of compact sets.
By choosing different collections of sets of similar type, one may use them to produce very general fractals. 
The collections of compact sets that produce self-similar fractals are simple to generate, whereas the collections that produce non-self-similar fractals are often more complicated. 
In this paper, we introduce the notion of an $\mathcal{F}$-limit set. 
We propose that it provides a simple, systematic way to represent collections of compact sets that produce general fractals including selfsimilar and non-selfsimilar ones, as well as deterministic and random ones. 

When $X$ is a metric space, a limit set is defined for a certain collection of subsets of $X$ indexed by the nodes of a tree as described as follows.
Let $\{n_k\}_{k = 1}^\infty$ be a sequence of positive integers. Let $D_0 = \emptyset$ and for each $k\ge 1$, let
\begin{equation}
D_k := \{(i_1, \dots, i_k) : 1 \leq i_j \leq n_k, 1 \leq j \leq k\}
\end{equation}
be the collection of all {\it words} of length $k$ with letters from the alphabet $\{1, \dots, n_k\}$.
With such a collection, we let $D := \cup_{k = 0}^\infty D_k$. 
The collection $D$ has a naturally directed tree structure (see Figure~\ref{tree 1}), where $k$ represents the generation, and $n_k$ denotes the number of children in generation $k$ that each parent set from generation $k-1$ has.
We call $D$ a {\it tree generated by $\{n_k\}_{k=1}^\infty$} or simply, a {\it tree}. 
If $n_k = m$ for all $k \in \N$ we say that $D$ is an {\it $m$-ary tree}.
\begin{figure}[h]
\begin{tikzpicture}[level distance=1.2cm,sibling distance=.001cm, 
   edge from parent/.style={draw,-latex}]
\Tree[.$\emptyset$ 
        [.$(1)$ 
            [.$(1,1)$ 
                $(1,1,1)$ $(1,1,2)$
            ]
            [.$(1,2)$
                $(1,2,1)$ $(1,2,2)$
            ]
        ]
        [.$(2)$ 
            [.$(2,1)$ 
                $(2,1,1)$ $(2,1,2)$
            ]
            [.$(2,2)$
                $(2,2,1)$ $(2,2,2)$
            ]
        ]
    ]
\end{tikzpicture}
\caption{Generations $k = 0, 1, 2$, and $3$ of the tree structure of a $2$-ary tree $D$.}
\label{tree 1}
\end{figure}
For any two words $\sigma = (\sigma_1, \dots, \sigma_k) \in D_k$ and $\tau = (\tau_1, \dots, \tau_i) \in D_i$, we define
\begin{equation}\label{dksigma}
\sigma \ast \tau := (\sigma_1, \dots, \sigma_k, \tau_1, \dots, \tau_i)\in D_{k+i}.
\end{equation}

\begin{definition}\label{def: limit set}
Given a tree $D$, and a collection $\mathcal{J} := \{J_\sigma : \sigma \in D\}$ of subsets of a metric space $X$, the {\it limit set} $F$ of $\mathcal{J}$ is defined to be 
\begin{equation}\label{eqn: Limit_set}
F := \bigcap_{k \geq 0} E_k \quad \text{ where } \quad E_k := \bigcup_{\sigma \in D_k} J_\sigma. 
\end{equation}
\end{definition}

When producing general fractals using limit sets, one typically chooses the elements $J_\sigma$ in the collection $\mathcal{J}$ to be alike. 
For instance, for the well-known Moran sets, every $J_\sigma$ is required to be ``similar'' to the root $J_\emptyset$.
Since their introduction by Moran~\cite{Moran}, Moran sets have been studied extensively by many authors with various approaches \cite[and references therein]{Beardon, inhom, Su, Li, Wen}.  
We reproduce the definition here with a more current interpretation.
\begin{definition}[\cite{Su}]\label{moran set}
Suppose that $J \subseteq\R^N$ is a compact set with nonempty interior.  
Let $\{n_k\}_{k \geq 1}$ be a sequence of positive integers, and $\{ \Phi_k\}_{k \geq 1}$ be a sequence of positive real vectors with 
\begin{equation}
\Phi_k = (c_{k,1}, c_{k,2}, \dots, c_{k,n_k}), \sum_{1 \leq j \leq n_k} c_{k,j} \leq 1, k \in \N.
\end{equation} 
Suppose that $\mathcal{J} := \{J_\sigma : \sigma \in D\}$ is a collection of subsets of $\R^N,$ where $D$ is the tree generated by $\{n_k\}_{k = 1}^\infty$.  
 We say that the collection $\mathcal{J}$ fulfills the Moran Structure provided it satisfies the following Moran Structure Conditions (MSC):
\begin{enumerate}
\item[MSC(1)] $J_{\emptyset} = J.$
\item[MSC(2)] For any $\sigma \in D$, $J_\sigma$ is geometrically similar to $J$.  That is, there exists a similarity $S_\sigma: \R^N \to \R^N$ such that $J_\sigma = S_\sigma(J)$. 
\item[\quad MSC(3)] For any $k \geq 0$ and $\sigma \in D_k$, $J_{\sigma \ast 1}, \dots , J_{\sigma \ast n_k}$ are subsets of $J_\sigma$, and $int(J_{\sigma \ast i}) \cap int(J_{\sigma \ast j}) = \emptyset$ for $i \neq j$.
\item[MSC(4)] For any $k \geq 1$ and $\sigma \in D_{k-1}, 1 \leq j \leq n_k,$ 
\begin{equation}
\label{ckj}
\dfrac{diam(J_{\sigma \ast j})}{diam(J_{\sigma})} = c_{k,j}.	
\end{equation}
\end{enumerate}
If $\mathcal{J}$ satisfies the Moran Structure Conditions, then we call the limit set of $\mathcal{J}$ a {\it Moran set}.
\end{definition}

Using the limit sets of the collections $\mathcal{J} := \{J_\sigma : \sigma \in D\}$
that satisfy the MSC is a great way to produce geometrically similar fractals (including self-similar and generalized self-similar fractals~\cite{Su}). 
A popular way to generate such collections is by using iterated function systems (IFSs).
An {\it iterated function system} (IFS) on $X$ is a finite family $\{S_1, \dots, S_m\}$ of contractions on $X$, where a {\it contraction on $X$} is a Lipschitz function from $X$ to $X$ with Lipschitz constant strictly less than $1$ (see~\cite{Falconer} for more details and applications). 
A nonempty compact subset $F \subseteq X$ is called an {\it attractor} of an IFS $\{S_1, \dots, S_m\}$ if $F = \bigcup_{i = 1}^m S_i(F)$. 
For example, the $\frac{1}{3}$-Cantor set is the attractor of the IFS $\{x/3, (x + 2)/3\}$ on $[0, 1]$. 

As shown in ~\cite{hutch}, the attractor of an IFS $\{S_1, \dots, S_m\}$ on a compact metric space $X$ is given by the limit set of the collection $\mathcal{J} = \{J_\sigma : \sigma \in D\}$ of compact subsets of $X$ defined by 
\begin{equation}\label{eq: IFS_LimitSet}
J_\emptyset = X \quad \text{ and } \quad J_\sigma = S_{i_1}\circ \dots \circ S_{i_m}(J_\emptyset)
\end{equation}
for all $\sigma = (i_1, \dots, i_m) \in D$, where $D$ is an $m$-ary tree. 
If in addition, $S_1, \dots, S_m$ are similarities, then $\mathcal{J}$ satisfies the Moran Structure Conditions by letting $n_k = m$ and $S_\sigma := S_{i_1}\circ \dots \circ S_{i_m}$ for $\sigma = (i_1, \dots, i_m)$ in MSC(2). 
The resulting Moran set is self-similar and agrees with the attractor of the IFS $\{S_1, S_2, \cdots, S_m\}$.  
The dimension of the limit set can be quickly calculated from the Moran-Hutchinson formula in \cite{hutch}.  

IFSs provide simple procedures for constructing the collections $\mathcal{J}$ used to generate self-similar fractals. By construction, each element $J_\sigma$ is similar to the root-ancestor $J_\emptyset$.  Nevertheless, the fractals that we see in nature are not necessarily strictly self-similar. Within a fixed scale, the fractals may look different at different places.
In these fractals, each child $J_{\sigma*i}$ is typically similar to its parent $J_\sigma$ with a small variation.  
Accumulations of these variations after many generations can cause a larger variation between each element $J_\sigma$ with its root-ancestor $J_\emptyset$. 
A natural question is: how to mathematically model these general fractals? 
We are looking for a generating method that has the following attributes.
\begin{itemize}
    \item[(A)] The method generates general fractals including selfsimilar and non-selfsimilar ones, as well as deterministic and random ones; 
    \item[(B)] Any child in each generation is nearly similar to its parent, i.e., each child is similar, but with the possibility of a ``small" variation;
    \item[(C)] The method is systematic and computationally simple.
\end{itemize}
Note that, (B) means that each child is obtained from its parent under some nearly similar map. 
In other words, a method with attribute (B) preserves the tree structure of $D$.

The limit sets generated from some collections of subsets are good candidates that satisfy (A), while the IFS procedure is an ideal one that satisfies (C). 
Towards a method that will also satisfy (B), a natural attempt is to combine them; by using limit sets from the collections $\mathcal{J}$ given in (\ref{eq: IFS_LimitSet}) after some perturbations of the IFSs. 
However, by~(\ref{eq: IFS_LimitSet}), it follows that $J_{\sigma* \tau}=S_\sigma (J_\tau)$ for any $\sigma, \tau \in D$, and in particular,
\[J_{\sigma*i}=S_\sigma (S_i(J_\emptyset))\ne S_i (S_\sigma(J_\emptyset))=S_i (J_\sigma).\]
That is, the $i^{th}$ child $J_{\sigma*i}$ of the $J_{\sigma}$ in the tree structure $D$ is not generated by applying $S_i$ to parent $J_\sigma$.
In this sense, the IFS procedure does not preserve the tree structure $D$. 
To keep the tree structure, we would like to obtain $J_{\sigma \ast i}$ by applying some function $\tilde{S}_i$ on $J_\sigma$. 
In this case, we would need $\tilde{S}_i = S_{\sigma}\circ S_i \circ S_\sigma^{-1}$ since
\[
J_\sigma \xrightarrow{S_\sigma^{-1}} J_\emptyset \xrightarrow{S_{i}} J_i \xrightarrow{S_{\sigma}} J_{\sigma\ast i}.
\]
Here the child $J_{\sigma*i}$ is ``forced'' to be similar to the root-ancestor $J_\emptyset$ under $S_{\sigma}\circ S_i$, rather than depending only on its parent $J_\sigma$. 
This feature causes trouble when one tries to generate general fractals that have variations among different scales and places.  

To overcome this issue, in this article, we modify the above attempt by considering certain types of set-valued mappings (see Definition~\ref{compression operator}) that directly map $J_{\sigma}$ to $J_{\sigma*i}$. 
These set-valued mappings can be used to generate a type of limit set, called an $\mathcal{F}$-limit set.
An $\mathcal{F}$-limit set is determined by a sequence of points in a fixed parameterization space $M$.
By choosing different types of sequences in $M$, we are able to get
various types of fractals: from self-simlilar to non self-similar, and from deterministic to random.
Standard fractals that can be obtained by IFSs correspond to constant sequences in $M$, whereas non-constant sequences produce non-self-similar fractals. 
Additionally, since the general process of constructing an $\mathcal{F}$-limit set is independent of the sequence, the computational complexity of producing an $\mathcal{F}$-limit set from a general sequence is the same as that of an IFS obtained from a constant sequence. As a result, our $\mathcal{F}$-limit set approach satisfies all the attributes (A), (B) and (C) listed above.

Another novelty of the article is the estimation of the Hausdorff dimension of limit sets.
In \S\ref{sec: Hausdorff} we find bounds for the Hausdorff dimension of the limit sets in a general metric space setting of a collection of bounded sets, not necessarily satisfying the MSC conditions. 
In particular, we introduce the concept {\it uniform covering condition} in Definition \ref{def: ucc} for the purpose of studying the lower bound of the Hausdorff dimension of the limit set, and provide sufficient conditions for this condition in later sections.

The article is organized as follows. 
After studying the Hausdorff dimension of limit sets in \S\ref{sec: Hausdorff}, we systematically formulate the general setup for the construction of $\mathcal{F}$-limit sets in \S 3. 
The Hausdorff dimensions of $\mathcal{F}$-limit sets are then estimated in \S \ref{dimensions}.
After that, in \S \ref{sec: examples} we apply the results to specific examples, including modifications of the Cantor set, the Sierpi\'nski triangle, and the Menger sponge.  
We also give a remark to discuss similarities and differences of this construction with $V-$variable fractals created by Barnsley, Hutchinson, and Stenflo in \cite{vvar1}, \cite{vvar2}. 
In particular, our $\mathcal{F}$-limit sets are analogous to $\infty-$variable fractals.
In \S \ref{UCC}, we explore the sufficient conditions needed for a fractal to satisfy the uniform covering condition, which plays a vital role in computing a lower estimate for the Hausdorff dimension of a limit fractal.

\section{Hausdorff Dimension of the Limit Sets}\label{sec: Hausdorff}
In this section we investigate the Hausdorff dimension $\dim_H(F)$ of the limit set $F$ defined in (\ref{eqn: Limit_set}) of a collection $\mathcal{J}$ that does not necessarily satisfy all the MSC conditions.
To start, we determine an upper bound for the dimension of the limit set $F$ by considering the step-wise relative ratios between the diameters of sets.

\begin{proposition}\label{thm: upper}
	Suppose $\mathcal J := \{J_{\sigma} : \sigma \in D\}$ is a collection of bounded subsets of a metric space $(X,d)$, and $s>0$. 
	Let	$E_{k} = \bigcup_{\sigma \in D_{k}} J_{\sigma}, \text{ and } F = \bigcap_{k \geq 0} E_{k}$ be defined as in (\ref{eqn: Limit_set}).  If there exists a sequence of positive numbers $\{c_k\}_{k=1}^{\infty}$ such that
	\[\liminf_{k\rightarrow \infty} \prod_{i=1}^k c_i =0\]
	and
	\begin{equation}
	\label{eqn: c}
	\sum_{j=1}^{n_k} \left (\diam(J_{\sigma*j})\right)^s \leq c_k \left( \diam(J_\sigma)\right)^s,
	\end{equation}
	for all $\sigma\in D_{k-1}$ and all $k=1,2,\cdots$, then $dim_H(F) \le s$.
\end{proposition}

\begin{proof}
	We prove by using mathematical induction that for $k=1,2,\cdots,$
	\begin{equation}
	\label{general_form}
	\sum_{\sigma \in D_k}(\diam(J_{\sigma}))^s \le \left(\prod_{i=1}^k c_i \right) (diam(J_\emptyset))^s.
	\end{equation}
	When $k=1$, (\ref{general_form}) follows from (\ref{eqn: c}). Now assume (\ref{general_form}) is true for some $k\ge 1$.  Then by (\ref{dksigma}), (\ref{eqn: c}), and (\ref{general_form}),
	
	\begin{eqnarray*}
		\sum_{\sigma\in D_{k+1}} (\diam(J_{\sigma}))^s &=&\sum_{\sigma\in D_k} \left(\sum_{j=1}^{n_{k+1}} (\diam(J_{\sigma*j}))^s\right)\\
		&\le& c_{k+1} \sum_{\sigma\in D_k}(\diam(J_{\sigma}))^s \le\left(\prod_{i=1}^{k+1} c_i \right)(\diam(J_\emptyset))^s
	\end{eqnarray*}
as desired. By the induction principle, (\ref{general_form}) holds for all $k=1,2,\cdots.$
	For each $k$, set
	\[\delta_k=\max\{\diam(J_\sigma): \sigma\in D_k\}>0.\]
	Then, by (\ref{general_form}),  $\delta_k\le  \left(\prod_{i=1}^k c_i \right)^{1/s}\diam(J_\emptyset)$.  Moreover, by (\ref{general_form})
	\[\mathcal{H}_{\delta_k}^{s}(F) \le \mathcal{H}_{\delta_k}^{s}(E_k)\le \sum_{\sigma\in D_k} \alpha(s)\left(\frac{\diam(J_\sigma)}{2}\right)^s \le \left(\prod_{i=1}^k c_i \right)\alpha(s)\left(\frac{\diam(J_\emptyset)}{2}\right)^s.\]
	Since $\liminf_{k\rightarrow \infty} \prod_{i=1}^k c_i =0$, there exists a sequence $\{k_t\}_{t=1}^\infty$ such that
	\begin{equation}\lim_{t\rightarrow \infty} \prod_{i=1}^{k_t} c_i =0.\end{equation} Thus, $\lim_{t\rightarrow \infty}\delta_{k_t}= 0$,
	$\mathcal{H}^{s}(F) =\lim_{t\rightarrow \infty}\mathcal{H}_{\delta_{k_t}}^{s}(F)=0 $, and hence $\dim_H(F) \le s$.
\end{proof}
Conversely, to study the lower bound on the Hausdorff dimension of the limit set $F$, we introduce the following concept.

\begin{definition}[uniform covering condition]
\label{def: ucc}
		Let $\mathcal J := \{J_{\sigma} : \sigma \in D\}$ be a collection of compact subsets of a metric space $(X, d)$, and $F$ be the limit set of $\mathcal{J}$ as given in (\ref{eqn: Limit_set}).  $\mathcal J$ is said to satisfy the \textit{uniform covering condition} if there exists a real number $\gamma>0 $ and a natural number $N$ such that for any closed ball $B$ in $X$, there exists a subset $D_B\subseteq D$ with cardinality of $D_B$ at most $N$,
	\begin{equation}
	\label{condtion_lower}
	B\cap F \subseteq \bigcup_{\sigma\in D_B} J_\sigma  \text { and }  diam(B)\ge \gamma \sum_{\sigma\in D_B}diam(J_\sigma).
	\end{equation}
\end{definition} 

\begin{proposition}\label{thm:lower}
Let $\mathcal J := \{J_{\sigma} : \sigma \in D\}$ be a collection of compact subsets of a metric space $(X, d)$ with diam$(J_{\emptyset})>0$, and $F$ be the limit set of $\mathcal{J}$ as given in (\ref{eqn: Limit_set}).  If $\mathcal J$ satisfies the uniform covering condition, and if for some $s>0$,
	\begin{equation}
	\label{lower_eqn}\sum_{j=1}^{n_k} \diam(J_{\sigma*j})^s
	\geq  \diam(J_\sigma)^s
	\end{equation}
	for all $\sigma\in D_{k-1}$ and all $k=1,2,\cdots$, then $dim_H(F)\ge s$.
\end{proposition}

\begin{proof}
	
	We first show that under condition (\ref{lower_eqn}), there exists a probability measure $\mu$ on $X$ concentrated on $F$ such that for each $\sigma \in D$,
	\begin{equation}\label{mu_c} 
	\mu(J_\sigma)\le \left(\frac{\diam(J_\sigma)}{\diam(J_\emptyset)}\right)^s.
	\end{equation}	
	Let $\mu(J_\emptyset)=1$, and for each $\sigma \in D_k$ for $k>0$ and $i=1, \cdots, n_k$, we inductively set
	\[\mu(J_{\sigma *i}) =\frac{\diam(J_{\sigma*i})^s}{\sum_{j=1}^{n_k} \diam(J_{\sigma *j})^s} \mu(J_{\sigma}).\]
	For any Borel set $A$ in $X$, define
	\[ \mu(A) = \inf \left \{ \sum_{i=1}^{\infty} \mu(J_{\sigma_i}) : A \cap F \subseteq \bigcup_{i=1}^{\infty} J_{\sigma_i} \text{ and } J_{\sigma_i} \in \mathcal J \right \}. \]	
	One can check that $\mu$ defines a probability measure on $X$, concentrated on $F$.
	
	To prove (\ref{mu_c}) for $J_{\sigma}$, $ \forall \sigma \in D$, we proceed by using induction on $k$ when $\sigma \in D_k$. It is clear for $k=0$.  Now assume that (\ref{mu_c}) holds for each $\sigma \in D_k$ for some $k$.  Then by induction assumption and (\ref{lower_eqn}), for each $i=1,\cdots, n_{k+1}$,
	\begin{align*}
	\mu(J_{\sigma*i}) &=\frac{\diam(J_{\sigma*i})^s}{\sum_{j=1}^{n_k} \diam(J_{\sigma*j})^s} \mu(J_{\sigma}) \\
	&\leq \frac{\diam(J_{\sigma*i})^s}{\sum_{j=1}^{n_k} \diam(J_{\sigma*j})^s} \left(\frac{\diam(J_{\sigma})}{\diam(J_{\emptyset})}\right)^s \leq \left(\frac{\diam(J_{\sigma*i})}{\diam(J_\emptyset)}\right)^s.
	\end{align*}
	This proves inequality (\ref{mu_c}).
	
	Now, for any $\delta>0$, let $\{B_i\}$ be any collection of closed balls with $\diam(B_i)\le \delta$ and $F\subseteq \cup_i B_i$. For each $i$, let $D_{B_i}$ be the subset of $D$ corresponding to $B_i$ as given in equation (\ref{condtion_lower}). Note that
	\[ F\subseteq \bigcup_i B_i\cap F \subseteq \bigcup_i  \bigcup_{\sigma\in D_{B_i}}J_\sigma=\bigcup_{\sigma\in \tilde{D}}J_\sigma,
	\]		
	where $\tilde{D}:=\cup_{i=1}^{\infty} D_{B_i}\subseteq D$. 
	
	Let
	\begin{eqnarray*}
		C(s)&:=&\max\{\sum_{i=1}^N \left(x_i\right)^s: (x_1,x_2,\cdots, x_N)\in [0,1]^N \text{ with } \sum_{i=1}^N x_i=1\}\\
		&=&
		\begin{cases}
			N^{1-s}, &\text{ if }0<s<1 \\
			1, &\text{ if }s\ge 1.
		\end{cases}
	\end{eqnarray*}
	and $c(s)=\frac{\alpha(s)}{C(s)}\left(\frac{\gamma\diam(J_\emptyset)}{2}\right)^s>0$. Then, by (\ref{condtion_lower}) and (\ref{mu_c}),
\comments{	
	\begin{eqnarray*}
		&&\sum_i  \alpha(s)\left(\frac{\diam(B_i)}{2}\right)^s \\ &\ge & \sum_i  \frac{\alpha(s)}{2^s} \left(\gamma \sum_{\sigma\in D_{B_i}}\diam(J_\sigma)\right)^s \\
		&\ge &  \sum_i  \frac{\alpha(s)}{2^s C(s)}\gamma^s\sum_{\sigma\in D_{B_i}}\left(\diam(J_\sigma)\right)^s \\
		&\ge &  \frac{\alpha(s)}{2^s C(s)}\gamma^s\sum_{\sigma\in \tilde{D}}\left(\diam(J_\sigma)\right)^s \\
		&\ge &  \frac{\alpha(s)}{2^s C(s)}\gamma^s\left(\diam(J_\emptyset)\right)^s\sum_{\sigma\in \tilde{D}}\mu(J_\sigma) \\
		&\ge & c(s)\mu\left(\sum_{\sigma\in \tilde{D}}J_\sigma\right) \ge c(s)\mu(F) = c(s).
	\end{eqnarray*}
}
 	\begin{eqnarray*}
		&&\sum_i  \alpha(s)\left(\frac{\diam(B_i)}{2}\right)^s \ge  \sum_i  \frac{\alpha(s)}{2^s} \left(\gamma \sum_{\sigma\in D_{B_i}}\diam(J_\sigma)\right)^s \\
		&\ge &  \sum_i  \frac{\alpha(s)}{2^s C(s)}\gamma^s\sum_{\sigma\in D_{B_i}}\left(\diam(J_\sigma)\right)^s \ge   \frac{\alpha(s)}{2^s C(s)}\gamma^s\sum_{\sigma\in \tilde{D}}\left(\diam(J_\sigma)\right)^s \\
		&\ge &  \frac{\alpha(s)}{2^s C(s)}\gamma^s\left(\diam(J_\emptyset)\right)^s\sum_{\sigma\in \tilde{D}}\mu(J_\sigma) \ge  c(s)\mu\left(\sum_{\sigma\in \tilde{D}}J_\sigma\right) \ge c(s)\mu(F) = c(s).
	\end{eqnarray*}
	
	Thus, $\mathcal{H}^{s}(F) =\lim_{\delta \rightarrow 0}\mathcal{H}_{\delta}^{s}(F)\ge c(s)>0 $, and hence $\dim_H(F) \ge s$.
\end{proof}

Next, we provide a sufficient condition for $\mathcal J$ to satisfy the uniform covering condition in a general situation. 
Later in \S\ref{UCC}, we will provide another one with Theorem~\ref{thm: sufficient_uniform_covering} that is particularly useful for the examples we will encounter in \S\ref{sec: examples}.

Recall that a metric space $(X,d)$ is called doubling if there is some doubling constant $M>0$ such that any ball $B(x,r)$ in $X$ can be covered by at most $M$ balls $B(x_i, r/2)$ in $X$. Equivalently \cite[Lemma 2.3]{Hytonen}, $(X,d)$ is doubling if for any $\epsilon>0$, there exists a natural number $N_\epsilon$ such that for any $\rho>0$, any ball in $X$ of diameter $\rho$ contains at most $N_\epsilon$ many disjoint balls of diameter $\epsilon \rho$. Clearly, any Euclidean space is a doubling metric space.

 \begin{proposition}
	
	Let $\mathcal J := \{J_{\sigma} : \sigma \in D\}$ be a collection of compact subsets of a doubling metric space $(X, d)$, and $F$ be the limit set of $\mathcal{J}$ as given in (\ref{eqn: Limit_set}).  Suppose that $\mathcal J$ satisfies the following conditions: 
	\begin{enumerate}
		\item	 there exists a number $r\in(0,1]$ such that for any $k\in \mathbb{N}$ and for each $\sigma\in D_k$,
		\[rc_k\le diam(J_\sigma)\le \frac{c_k}{r}\]
		where $c_k:=\min\{	diam(J_{\bar{\sigma}}): \bar{\sigma}\in D_{k-1}\}$.
		\item there exists a number $\tau\in (0, 1]$ such that for each $\sigma\in D$, the convex hull of $J_\sigma$ contains a closed ball $W_\sigma$
		such that 
		\[diam(W_\sigma)\ge \tau \cdot diam(J_\sigma)\]
		and for each $k\in \mathbb{N}$, the collection $\{W_\sigma: \sigma\in D_k\}$ are pairwise disjoint.
	\end{enumerate}
 Then $F$ satisfies the uniform covering condition (\ref{condtion_lower}).
\end{proposition}
\begin{proof}	
	  For any closed ball $B$ in $X$, let $k$ be the number such that
	\[	 c_{k+1}\le diam(B) <c_k\]
	where by convention, we set $c_0=\infty$. Let
	\[D_B:=\{\sigma\in D_k: B\cap F\cap J_\sigma \neq \emptyset\}.\]	
	Note that
	\[B\cap F = B \cap F \cap \bigcup_{\sigma\in D_k} J_\sigma \subseteq\bigcup_{\sigma\in D_B} J_\sigma.\]	Also for any $\sigma\in D_B$, since $diam(J_\sigma)\le \frac{c_k}{r}$ and $B\cap J_\sigma \neq \emptyset$, it follows that $J_\sigma\subseteq \bar{B}(x_0, \frac{r+2}{2r} c_k)$, where $x_0\in X$ is the center of the ball $B$.  Thus, $W_\sigma \subseteq \bar{B}(x_0, \frac{r+2}{2r} c_k)$. Let $\rho=\frac{r+2}{r}c_k$ and $\epsilon =\frac{r^2 }{r+2}\tau$, then
	\[diam(W_\sigma)\ge \tau\cdot diam(J_\sigma) \ge \tau r c_k = \epsilon \rho.\]
	Since $\{W_\sigma: \sigma\in D_B\}$ are pairwise disjoint and $(X,d)$ is doubling, the cardinality of $D_B$ is at most $N:=N_\epsilon$.  
	On the other hand, for $\gamma =\frac{r^2}{N}$,	it holds that
	\begin{equation}
	diam(B)\ge c_{k+1}\ge rc_k= \gamma N \frac{c_{k}}{r} \ge \gamma \sum_{\sigma\in D_B}\frac{c_{k}}{r} \ge \gamma \sum_{\sigma\in D_B}diam(J_\sigma).
	\end{equation}
	As a result, $\mathcal J$ satisfies the condition (\ref{condtion_lower}) as desired.
\end{proof}	

\section{General Setup of $\F$-Limit sets}
\label{general setup}
We now formalize the construction of general fractals using $\F$-limit sets. 

Let $\mathcal{A} \subseteq 2^X$ be a collection of subsets of $X$, let $\mathcal{O}(\mathcal{A})$ denote the collection of all operators that map $\mathcal{A}$ to $\mathcal{A}$, and let $D$ be an $m$-ary tree. 
Given a set $E_0 \in \mathcal{A}$ and an $\mathcal{O}(\mathcal{A})^m$-valued map on $D$
\begin{equation}\label{general-operator-valued-map}
f : \sigma \mapsto f_\sigma = (f_\sigma^{(1)}, \dots, f_{\sigma}^{(m)}),
\end{equation}
we can construct the limit set of the collection $\mathcal{J}(f, E_0) := \{J_\sigma : \sigma \in D\}$ of subsets defined by
\begin{equation}\label{vertex map}
J_\emptyset := E_0, \quad \text{ and } \quad J_{\sigma\ast j} := f_\sigma^{(j)}(J_\sigma)
\end{equation}
for all $\sigma \in D$ and all $1 \leq j \leq m$. 
The corresponding set 
\begin{equation}\label{def: limit-set}
F = \bigcap_{k \geq 0} E_k \quad \text{ where } \quad E_k := \bigcup_{\sigma \in D_k} J_\sigma
\end{equation}
is called the {\it limit set generated by $f$ with initial set $E_0$}. 

\begin{remark}
    To ensure that our limit sets are not empty, the class $\mathcal{O}(\mathcal{A})$ of operators that we will be using is the class $\mathcal{C}(\mathcal{X})$ of compression operators defined as follows.
    For more general maps, instead of considering the intersection of $\{E_k\}$, one may study the limit of the sequence $\{E_k\}$ in some suitable metric.
    We leave this path of exploration to future research. 
\end{remark}

\begin{definition}\label{compression operator}
Given a collection $\mathcal{X}$ of compact subsets of $X$, an operator $f : \mathcal{X} \to \mathcal{X}$ is called a {\it compression operator} on $\mathcal{X}$ if $f(E) \subseteq E$ for all $E \in \mathcal{X}$.
We let $\mathcal{C}(\mathcal{X})$ denote the collection of all compression operators on $\mathcal{X}$.
\end{definition}

\begin{example}\label{ex:contractions}
Here are some examples of simple collections of compression operators. 
\begin{itemize}
\item[(a)] Let $S:X\rightarrow X$ be a contraction map on the metric space $X$. Let $\mathcal{X}_{S}$ be the collection of all compact subsets of $(X,d)$ with $S(E)\subseteq E$. Then $S$ is a compression on $\mathcal{X}_S$.

\item[(b)] Let $\mathcal{X}_{cpt}$ be the collection of all compact subsets of $(X,d)$, and let $K$ be a compact set. Then, $f_K(E):=E\cap K$ defines a compression $f_K$ on $\mathcal{X}_{cpt}$.
Note that $f_K$ is usually not given by a contraction map. 

\item[(c)] Let  $\mathcal{X}_I = \{ [a,b] : a,b \in \R\}$ be the collection of all closed intervals in $\R$. Then for each $p\in [0,1]$, both $f^{(1)}_p([a,b]):=[a, p(b-a)+a]$ and $f^{(2)}_p([a,b]):=[p(b-a)+a, b]$ define compressions $f^{(1)}_p$ and $f^{(2)}_p$ on $\mathcal{X}_I$. 

\item[(d)] Let $f$ be a compression on $\mathcal{X}$, and $K\in \mathcal{X}$. Then the ``restriction'' $f\vert_{K}$ of $f$ on $K$ defines a compression on
$\mathcal{X}\vert_K:=\{E\cap K: E\in \mathcal{X} \}$ because $f(E\cap K)\subseteq f(E)\cap f(K)\subseteq E\cap K$.
\end{itemize}
\end{example}
In this article, we study the $\mathcal{C}(\mathcal{X})^m$-valued maps on $D$ that are defined by the composition of two maps
\begin{equation}\label{vertex map diagram}
\begin{tikzcd}
  D \arrow[r, dashrightarrow] \arrow[d, "\k"] & \mathcal{C}(\mathcal{X})^m \\
  M  \arrow[ru, "\mathcal{F}"] 
\end{tikzcd}
\end{equation}
that factors through a set $M$ that is used as a way to parameterize some subset of compression operators. 
Throughout this article, our parameterization space $M$ will often be some Cartesian product of the closed unit interval $[0, 1]$. 
A limit set generated by such a map $f = \mathcal{F}\circ \k$ with an initial set $E_0$ is called an {\it $\mathcal{F}$-limit set}.
Let us formally state the definition of an $\mathcal{F}$-limit set.

\comments{
\blue{\begin{definition}
Let $\mathcal{X}$ be a collection of compact subsets of $X$, and let $\mathcal{C}(\mathcal{X})$ denote the collection of all compressions on $\mathcal{X}$. For any $m\in \mathbb{N}$, a marking of $\mathcal{C}(\mathcal{X})^m$ is a map 
\[
\mathcal{F} : M \to \mathcal{C}(\mathcal{X})^m; \qquad p \mapsto (f_{p}^{(1)}, \dots, f_{p}^{(m)}),
\]
where $M$ is a nonempty set, called a {\it parameterization space}. 
\end{definition}

\begin{definition}\label{F limit set}
Let $\mathcal{F}:  M \to \mathcal{C}(\mathcal{X})^m$ be a marking of $\mathcal{C}(\mathcal{X})^m$, and $D$ be an $m$-ary tree. We say that a subset $F$ of $X$ is an $\mathcal{F}$-limit set if it is the limit set generated by the composition map $f=\mathcal{F}\circ \k$ for some map map $\k : D \to M$ with an initial set $E_0\in \mathcal{X}$. Namely, $F$ is the limit set of a collection $\{J_\sigma : \sigma \in D\}$ of compact subsets of $X$ defined by 
\begin{equation}\label{F-limit-set}
J_\emptyset := E_0, \quad \quad J_{\sigma \ast j} = f_{\k(\sigma)}^{(j)}(J_\sigma) \quad \text{ for all } \sigma \in D \text{ and all } 1 \leq j \leq m.
\end{equation}

\end{definition}
}
}
\begin{definition}\label{F limit set}
Let $D$ be an $m$-ary tree, let $\mathcal{X}$ be a collection of compact subsets of $X$, and let $\mathcal{C}(\mathcal{X})$ denote the collection of all compressions on $\mathcal{X}$. 
We say that a subset $F$ of $X$ is an $\mathcal{F}$-limit set if there is a set $M$ and maps $\k : D \to M$, and 
\[
\mathcal{F} : M \to \mathcal{C}(\mathcal{X})^m; \qquad p \mapsto (f_{p}^{(1)}, \dots, f_{p}^{(m)})
\]
such that $F$ is the limit set of the collection $\{J_\sigma : \sigma \in D\}$ of compact subsets of $X$ defined by 
\begin{equation}\label{F-limit-set}
J_\emptyset := E_0, \quad \quad J_{\sigma \ast j} = f_{\k(\sigma)}^{(j)}(J_\sigma) \quad \text{ for all } \sigma \in D \text{ and all } 1 \leq j \leq m.
\end{equation}
The set $M$ is called a {\it parameterization space} and $\mathcal{F}$ is called a marking of $\mathcal{C}(\mathcal{X})^m$.
\end{definition}

In this article, we are interested in fixing a marking $\mathcal{F}$ and investigating the different $\mathcal{F}$-limits sets that can be generated by varying $\k$.
When $\mathcal{F}$ is fixed, we say that the limit set $F$ of Definition~\ref{F limit set} is an 

\begin{center}
{\it $\mathcal{F}$-limit set generated by the map $\k$ with initial set $E_0$.}
\end{center}
The following is an example of a marking that will produce different Cantor-like fractals when choosing different $\k$. 
\begin{example}\label{ex: cantor 1}
Let $M = [0, 1]^2$ and let $\mathcal{X}_I := \{[a, b] : a, b \in \mathbb{R}\}$. 
Define a marking $\mathcal{F}$ of $M$ into $\mathcal{C}(\mathcal{X}_I)^2$ by 
\[
\F(p) = (f^{(1)}_p, f^{(2)}_p) \quad \text{ for all } \quad p = (p_1, p_2) \in [0, 1]^2
\]
where 
\[
	f^{(1)}_p([a,b]) := [a, p_1(b-a)+a] \quad \text{ and } \quad f^{(2)}_p([a,b]) := [p_2(a-b)+b, b].
\]
By choosing $\k : D \to M$ to be the constant map defined by
\[
\k(\sigma) = (1/3, 1/3)\quad \text{ for all } \quad \sigma \in D,
\]
we get that the $\mathcal{F}$-limit set generated by $\mathcal{F}\circ \k$ with the initial set $E_0 = [0, 1]$ is equal to the standard $\frac{1}{3}$-Cantor set.
In Example~\ref{example: 1}, we give an example of a non-constant function $\k : D \to M$ that will produce an $\mathcal{F}$-limit set that has the same basic shape as the standard $\frac{1}{3}$-Cantor set, but is not self-similar. 
\end{example}

\begin{example}\label{ex: IFS}
For any IFS $S = \{S_1, \dots, S_m\}$ on a metric space $X$, let $\mathcal{J}_S := \{J_\sigma : \sigma \in D\}$, where each $J_\sigma$ is defined as in~(\ref{eq: IFS_LimitSet}) and $D$ is the $m$-ary tree.
Suppose there exists one collection $\mathcal{X}$ containing $\mathcal{J}_S$, and a marking 
\[
\F : M \to \mathcal{C}(\mathcal{X})^m; \quad p \mapsto (f_p^{(1)}, \dots, f_p^{(m)})
\]
from a parameterization space $M$ to $\mathcal{C}(\mathcal{X})^m$ such that for some $p^*\in M$,
\[
f_{p^*}^{(j)}(J_\sigma) = J_{\sigma\ast j}
\]
for all $\sigma\in D$ and $j = 1, \dots, m$.
Then, the $\mathcal{F}$-limit set generated by the constant map $\k(\sigma) := p^*$ agrees with the attractor of the IFS.
We will provide explicit examples of these types of markings in \S5. 
\end{example}
\comments{
\begin{example}\label{ex: IFS}
Let $S = \{S_1, \dots, S_m\}$ be an IFS on a compact metric space $X$ that satisfies the open set condition with $X$, and let $D$ be the $m$-ary tree. 
Define $\mathcal{J}_S = \{J_\sigma : \sigma \in D\}$, where the $J_\sigma$ are defined as in~(\ref{eq: IFS_LimitSet}).
First, for each $j = 1, \dots, m$ define the compression $f^{(j)} : \mathcal{J}_S \to \mathcal{J}_S$ by
\[
f^{(j)}(J_\sigma) := J_{\sigma\ast j} \quad \text{ for all } \sigma\in D.
\]
Notice that since $S$ satisfies the open set condition with $X$, by Lemma~\ref{IFS with open set condition}, the $f^{(j)}$ is well-defined for each $j = 1, \dots, m$. 
Now let $M$ be any non-empty set and fix some point $p_0 \in M$.
Choose any marking $\mathcal{F} : M \to \mathcal{C}(\mathcal{J}_S)^m$ of $M$ into $\mathcal{C}(\mathcal{J}_S)^m$ such that 
\[
\mathcal{F}(p_0) = (f^{(1)}, \dots, f^{(m)}).
\]
By choosing $\k$ to be the constant map $\k(\sigma) := p_0$ for all $\sigma \in D$, the $\mathcal{F}$-limit set generated by $\mathcal{F}\circ \k$ agrees with the attractor of the IFS.  
\end{example}
 \blue{*****delete above materials******}
 }

Since the $m$-ary tree $D$ is an ordered set, it is sometimes more convenient to represent the mapping $\k : D \to M$ as a sequence $\{\k_{\ell}\}_{\ell=0}^{\infty}$ in $M$ where $\k_\ell := \k(\sigma_\ell)$ for some ordering $\{\sigma_\ell\}_{\ell = 0}^\infty$ of $D$.  
This ordering of $D$ is given as follows: Set $\sigma_0 = \emptyset$ and for any $N\in \mathbb{N}$, define $\sigma_N := (i_1, \dots, i_k)$ where $1 \leq i_1, \dots, i_k \leq m$ are uniquely defined by the expression
\begin{equation}\label{ordering of D}
N=\sum_{p=0}^{k-1} m^p i_{k-p}. 
\end{equation}
In other words, the ordering of $D$ is determined by the map $\ell : D \to \mathbb{N}\cup\{0\}$ defined by $\ell(\emptyset)=0$ and 
\begin{equation}\label{ell_sigma} \ell(\sigma) = \sum_{p=0}^{k-1} m^p i_{k-p} \quad \text{ for all } \sigma=(i_1, i_2, \dots, i_k) \in D.
\end{equation}
Using this notation, we can rewrite Definition \ref{F limit set} as follows. 
To do so, let us introduce another notation. For each $n \ge 1$, let $\im(n)$ be the cardinality of $\cup_{k = 1}^n D_k$, namely
	\begin{equation} \label{eqn:p_m}
	\im(n)=m+m^2+\cdots+m^n=\frac{m^{n+1}-m}{m-1}.
	\end{equation}  
 Also set $\im(0)=0$.
 

\begin{definition}[Revision of Definition  \ref{F limit set}]
\label{def: revised}
Let $\mathcal{X}$ be a collection of compact subsets of $X$ and let $\mathcal{C}(\mathcal{X})$ denote the collection of all compressions on $\mathcal{X}$.  
Fix a map $\mathcal{F}$ from a non-empty set $M$ to $\mathcal{C}(\mathcal{X})^m$ denoted by $\mathcal{F}(p) = (f_{p}^{(1)}, \dots, f_{p}^{(m)})$.

For any sequence $\{\k_{\ell}\}_{\ell=0}^{\infty}$ in $M$ and $E_0 \in \mathcal{X}$, we iteratively define the sets   
\[
E_{m\ell+j} := f_{\k_\ell}^{(j)}(E_\ell) \in \mathcal{X}, \text{ for }\ell=0,1,2,\cdots,\ j=1,2,\cdots.
\]

	The limit set \begin{equation} 
	F=\bigcap_{n=1}^{\infty} \bigcup_{\ell=\im(n-1)+1}^{\im(n)} E_{\ell} \label{F_limit_set} \end{equation}
	is called the $\F$-limit set generated by the sequence $\{\k_\ell\}_{\ell = 0}^\infty$ with initial set $E_0$, where $\im(n)$ is defined as in~(\ref{eqn:p_m}). 
\end{definition}

Note that the $\mathcal{F}$-limit set generated by the sequence $\{\k_\ell\}_{\ell = 0}^\infty$ with initial set $E_0$ is the $\mathcal{F}$-limit set generated by $\mathcal{F}\circ \k$ where the map $\k : D \to M$ is given by $\sigma \mapsto \k_{\ell(\sigma)}$ and $\ell(\sigma)$ as defined in~(\ref{ell_sigma}).

\section{Hausdorff dimensions of $\F$-Limit sets}
\label{dimensions}

In this section, we fix an $m$-ary tree $D$, a parameterization space $M$, a collection $\mathcal{X}$ of compact subsets of a metric space $X$, and a marking $\mathcal{F}$ of $M$ into $\mathcal{C}(\mathcal{X})^m$. 

As indicated in Propositions \ref{thm: upper} and \ref{thm:lower}, the relative ratio between the diameters of the sets plays an important role in the calculation of the dimension of the limit set.  Therefore, we introduce the following definition.

\begin{definition}
	For any compression $g: \mathcal{X}\to \mathcal{X}$, define
	\begin{equation}
	U(g)=\sup_{E\in \mathcal{X}} \frac{\diam(g(E))}{\diam(E)}, \text{ and } L(g)=\inf_{E\in \mathcal{X}} \frac{\diam(g(E))}{\diam(E)}.
	\end{equation}
\end{definition}
Note that, for each $E\in \mathcal{X}$,
\begin{equation}
\label{LgUg}
L(g)\cdot \diam(E)\le \diam(g(E)) \le U(g) \cdot \diam(E).
\end{equation}

For any $p \in M$ and $f_p = (f_{p}^{(1)}, \cdots , f_{p}^{(m)}) \in \mathcal{C}(\mathcal{X})^m$, define
\[\mathbf{U}_p=\left( U(f_{p}^{(1)}), \cdots,  U(f_{p}^{(m)})\right) \in \mathbb{R}^m,\]
and 
\[ \mathbf{L}_p=\left( L(f_{p}^{(1)}), \cdots,  L(f_{p}^{(m)})\right) \in \mathbb{R}^m.\]

Also, for each $x=(x_1, \cdots, x_m)\in \mathbb{R}^m$ and $s>0$, denote
\[||x||_s=\left(\sum_{i=1}^m |x_i|^s\right)^{\frac{1}{s}}.\]

These notations, Proposition \ref{thm: upper} and Proposition \ref{thm:lower} motivate our main theorem.
\begin{theorem} \label{thm:ratio_bounds}
	Let $F$ be the $\F$-limit set generated by the map $\k : D \to M$ with initial set $J_{\emptyset}$, and $s> 0$. 
	\begin{enumerate}
		\item[(a)] If $F$ satisfies the uniform covering condition (\ref{condtion_lower}) and
		\[\inf_{\sigma \in D} \{||\mathbf{L}_{\k(\sigma)}||_{s}\}\ge 1,\]
		then $\dim_H(F) \ge s$.
		\item[(b)] If
		\[\sup_{\sigma \in D}\{||\mathbf{U}_{\k(\sigma)}||_{s}\}<1,\]
		then $\dim_H(F)\le s$.
	\end{enumerate}
\end{theorem}

\begin{proof}
	(a)  By (\ref{F-limit-set}) and (\ref{LgUg}), for all $\sigma \in D$, 
	\begin{eqnarray*}
		& &\sum_{j=1}^{m} \diam(J_{\sigma * j})^s =\sum_{j=1}^{m} \diam\left(f_{\k(\sigma)}^{(j)}(J_{\sigma})\right)^s 
		\ge  \sum_{j=1}^{m}\left(L(f_{\k(\sigma)}^{(j)})\right)^s \diam(J_{\sigma})^s \ge \diam(J_{\sigma})^s.
	\end{eqnarray*}
	Thus, by Proposition \ref{thm:lower}, $\dim_H(F)\ge s$.
	
	(b) Similarly, for all $\sigma \in D$, 
	\[\sum_{j=1}^{m} \diam(J_{\sigma*j})^s \le \sum_{j=1}^{m}\left(U(f_{\k(\sigma)}^{(j)})\right)^s \diam(J_{\sigma})^s \le c \cdot \diam(J_{\sigma})^s,\]
	where
	\[c:=\sup_\sigma\{(||\mathbf{U}_{\k(\sigma)}||_{s})^s\}<1.\]
	By Proposition \ref{thm: upper}, $\dim_H(F)\le s$.
\end{proof}
	 
In the following, we will use the notation from Definition~\ref{def: revised} to describe the construction of the $\F$-limit sets.  
Clearly, using this notation, Theorem~\ref{thm:ratio_bounds} simply says that if $F$ satisfies the uniform covering condition (\ref{condtion_lower}) and 	$\displaystyle \inf_{\ell} \{||\mathbf{L}_{\k_{\ell}}||_{s}\}\ge 1,$ then $\dim_H(F) \ge s$, and if $ \displaystyle \sup_{\ell}\{||\mathbf{U}_{\k_\ell}||_{s}\}<1$, then $\dim_H(F)\le s$.	

When both $\{||\mathbf{L}_{\k_\ell} ||_s\}_{\ell=0}^{\infty}$ and $\{||\mathbf{U}_{\k_\ell} ||_s\}_{\ell =0}^{\infty}$ are convergent sequences, the following corollary enables us to quickly estimate the dimension of $F$.

\begin{corollary}\label{cor:limit_sup_inf}
	Let $F$ be the $\mathcal{F}$-limit set generated by the sequence $\{ \k_{\ell} \}_{\ell =0}^{\infty}$ with initial set $E_0$.
	
		\begin{enumerate}
		\item[(a)] Let $\underline{s_*} := \sup \left\{s: \liminf_{\ell \rightarrow \infty}\{||\mathbf{L}_{\k_\ell} ||_s\}>1 \right\}.$  Then
		\begin{equation}
		\dim_H(F) \ge \underline{s_*},
		\end{equation}
		provided $F$ satisfies the uniform covering condition (\ref{condtion_lower}).
		
		\item[(b)] Let $\overbar{s^*} := \inf \left \{s: \limsup_{\ell \rightarrow \infty}\{||\mathbf{U}_{\k_\ell}||_s\}<1 \right\}$.  Then	\begin{equation}\label{eqn: s_star}
		\dim_H(F)\le \overbar{s^*}.
		\end{equation}
		
		\end{enumerate}

\end{corollary}
\begin{proof}
	For any $0< s< \underline{s_*}$, by the definition of $\underline{s_*}$,
	\[\liminf_{\ell \rightarrow \infty}\{||\mathbf{L}_{\k_\ell} ||_s\}>1.\]
	Thus, when $\ell_*\in \mathbb{N}$ is large enough,
	\[\inf_{\ell\ge \ell_*}\{||\mathbf{L}_{\k_\ell}||_s \} \ge 1, \qquad \text{i.e. } \inf_{\ell\ge 0}\{||\mathbf{L}_{\k_{\ell_*+\ell}}||_s\} \ge 1.\]
	Since $F\cap E_{\ell_*}$ is the set generated by
	the sequence $\{\k_{\ell_*+\ell} \}_{\ell =0}^{\infty}$ with initial set $E_{\ell_*}$, by Theorem \ref{thm:ratio_bounds}, it follows that $\dim_H(F\cap E_{\ell_*}) \ge s$ for any $\ell_*$ large enough. This implies that $\dim_H(F) \ge s$ for any $s<\underline{s_*}$ and hence $\dim_H(F) \ge \underline{s_*}$. Similarly, we also have $\dim_H(F)\le \overbar{s^*}$.
\end{proof}

In the following corollaries, we will see that bounds of the dimension of $F$ can also be obtained from corresponding bounds on $\mathbf{L}_{\k_\ell}$ and $\mathbf{U}_{\k_\ell}$. 

\textbf{Notation.} For any two points $x=(x_1,\cdots, x_m)$ and $y=(y_1,\cdots, y_m)$ in $\mathbb{R}^m$, we say $x \le y$ if $x_i \le y_i$ for each $i=1,\cdots, m$.

\begin{corollary} \label{cor: bound_on_L_U}
	Let $\mathbf{t}=(t_1,\cdots, t_m) $  and $\mathbf{r}=(r_1,\cdots, r_m) $ be two points in $(0,1)^m \subseteq \mathbb{R}^m$. Let $s_*$ and $s^*$ be the solutions to $||\mathbf{t}||_{s_*}=1$, and $||\mathbf{r}||_{s^*}=1$ respectively,  i.e. 
	\[t_1^{s_*}+t_2^{s_*}+\cdots+t_m^{s_*}=1,  \text{ and } r_1^{s^*}+r_2^{s^*}+\cdots+r_m^{s^*}=1.\]
 Let $F$ be the $\mathcal{F}$-limit set generated by the sequence $\{ \k_{\ell} \}_{\ell =0}^{\infty}$ with initial set $E_0$.
	\begin{enumerate} 
		\item[(a)] If $\mathbf{L}_{\k_\ell}\ge \mathbf{t}$ for all $\ell$ and $F$ satisfies the uniform covering condition (\ref{condtion_lower}), then $\dim_H(F)\ge s_*$.
		
		\item[(b)]If $\mathbf{U}_{\k_\ell} \le \mathbf{r}$ for all $\ell$, then
		$\dim_H(F) \le s^*$. 
		
		\item[(c)] If $ \mathbf{L}_{\k_\ell}=\mathbf{r} =\mathbf{U}_{\k_\ell}$ for all $\ell$ and $F$ satisfies the uniform covering condition (\ref{condtion_lower}), then $\dim_H(F)=s^*$.
	\end{enumerate}
\end{corollary}
\begin{proof}
	(a) Let $0< s<s_*$. Then,
	\[\inf_\ell\{||\mathbf{L}_{\k_\ell}||_{s}\}\ge ||\mathbf{t}||_s \ge ||\mathbf{t}||_{s_*}= 1.\]
	Thus, by Theorem \ref{thm:ratio_bounds}, $\dim_H(F)\ge s$ for any $s<s_*$, and hence $\dim_H(F)\ge s_*$.
	
	(b) Similarly,  let $0<s^*<s$. Then,
	\[\sup_\ell\{||\mathbf{U}_{\k_\ell}||_{s}\}\le ||\mathbf{r}||_s < ||\mathbf{r}||_{s^*}= 1.\]
	Thus, by Theorem \ref{thm:ratio_bounds},  $\dim_H(F)\le s$ for any $s>s^*$, and hence $\dim_H(F)\le s^*$.
	
	(c) follows from (a) and (b).
\end{proof}

A special case of Corollary \ref{cor: bound_on_L_U} gives the following explicit formulas for the bounds on the dimension of $F$.
\begin{corollary}\label{cor:uniform_bounds}
Let $F$ be the $\mathcal{F}$-limit set generated by the sequence $\{\k_{\ell}\}_{\ell =0}^{\infty}$ with initial set $E_0$.
Let
	\[\mathbf{t}=(t, \cdots, t) \text{ and } \mathbf{r}=(r, \cdots, r),\]
	for some $0<t, r<1$. 
	\begin{enumerate}
		\item[(a)] If $\mathbf{L}_{\k_\ell}\ge \mathbf{t}$ for all $\ell$ and $F$ satisfies the uniform covering condition (\ref{condtion_lower}), then $\dim_H(F)\ge \frac{\log{m}}{-\log t}$.
		
		\item[(b)]If $\mathbf{U}_{\k_\ell} \le \mathbf{r}$ for all $\ell$, then
		$\dim_H(F) \le \frac{\log{m}}{-\log r}$. 
		
		\item[(c)] If $ \mathbf{L}_{\k_\ell}=\mathbf{r} =\mathbf{U}_{\k_\ell}$ for all $\ell$ and $F$ satisfies the uniform covering condition (\ref{condtion_lower}), then $\dim_H(F)=\frac{\log{m}}{-\log r}$.
	\end{enumerate}
\end{corollary}

Other types of bounds on $\mathbf{L}_{\k_\ell}$ and $\mathbf{U}_{\k_\ell}$ can also be used to provide bounds on $\dim_H(F)$, as indicated by the following result.
\begin{corollary}\label{cor: mean}
	Let $F$ be the $\mathcal{F}$-limit set generated by the sequence $\{\k_{\ell} \}_{\ell =0}^{\infty}$ with initial set $E_0$. 
	\begin{enumerate}
		\item[(a)] If $F$ satisfies the uniform covering condition (\ref{condtion_lower}) and 
		\[w:=\inf_{\ell}\{||\mathbf{L}_{\k_\ell}||_1\} \ge 1,\]
		then $\dim_H(F) \ge \frac{\log(m)}{\log(m)-\log(w)}$.
		\item[(b)]
		If 
		\[u:=\sup_{\ell}\{||\mathbf{U}_{\k_\ell}||_1\} < 1,\]
		then $\dim_H(F) \le \frac{\log(m)}{\log(m)-\log(u)}$.
	\end{enumerate}
\end{corollary}
\begin{proof}
	(a). In this case, for $s = \frac{\log(m)}{\log(m)-\log(w)} \ge 1$, we have
	\[\frac{\sum_{j=1}^{m} \left(L\left(f_{\k_\ell}^{(j)}\right) \right)^{s}}{m}\ge \left(\frac{\sum_{j=1}^{m}L\left(f_{\k_\ell}^{(j)}\right)}{m}\right)^s \ge \left(\frac{w}{m}\right)^s\]
	for each $\ell$. Thus,
	\[\inf_\ell\{||\mathbf{L}_{\k_\ell}||_s\}\ge m^{\frac{1}{s}}\frac{w}{m} = 1,\]
	then by Theorem \ref{thm:ratio_bounds}, $\dim_H(F) \ge s$.
	
	(b).  In this case, for any $1 \ge s > \frac{\log(m)}{\log(m)-\log(u)}$, we have
	\[\frac{\sum_{j=1}^{m}\left(U\left(f_{\k_\ell}^{(j)}\right)\right)^{s}}{m}\le \left(\frac{\sum_{j=1}^{m}U\left(f_{\k_\ell}^{(j)}\right)}{m}\right)^s \le \left(\frac{u}{m}\right)^s\]
	for each $\ell$. Thus,
	\[\sup_\ell\{||\mathbf{U}_{\k_\ell}||_s\}\le m^{\frac{1}{s}}\frac{u}{m} <1.\]
	By Theorem \ref{thm:ratio_bounds}, $\dim_H(F) \le s$. Hence, $\dim_H(F) \le \frac{\log(m)}{\log(m)-\log(u)}$.
\end{proof}

Note that this corollary generally provides better bounds on $\dim_H(F)$ than those obtained from directly applying Theorem \ref{thm:ratio_bounds}.

\section{Examples of $\F$-Limit sets}\label{sec: examples}

In this section, we provide concrete examples of $\mathcal{F}$-limit sets in dimensions 1, 2, and 3.
Within each example, we provide a parameter space and a marking such that the $\mathcal{F}$-limit set generated by a certain constant sequence will result in the classical fractals: the Cantor set, the Sierpi\'nski triangle, and the Menger sponge. 
By choosing non constant sequences in our parameter spaces, we build non self-similar variations of the classical fractals mentioned above, and with the results of \S\ref{dimensions} and \S6, we are able to find (or estimate) their Hausdorff dimensions. 
It is worth reiterating that the computational complexity of the construction is independent of the sequence in the parameter space.

\begin{subsection}{Cantor-like sets}\label{cantor-like sets}
To construct our Cantor-like sets, in this subsection we choose $D$ to be the $2$-ary tree, $M := [0, 1]^2$, and 
\begin{equation}
\label{cantor_X}
    \mathcal{X} := \{[a, b] : a, b \in \mathbb{R}\}.
\end{equation}
Also, as in Example~\ref{ex: cantor 1}, define the marking $\mathcal{F}$ of $M$ into $\mathcal{C}(\mathcal{X})^2$ as 
\[
\F(p) = (f^{(1)}_p, f^{(2)}_p) \quad \text{ for all } \quad p = (p_1, p_2) \in [0, 1]^2
\]
where 
\[
	f^{(1)}_p([a,b]) := [a, p_1(b-a)+a] \quad \text{ and } \quad f^{(2)}_p([a,b]) := [p_2(a-b)+b, b].
\]

\comments{
	We first consider Cantor-like sets.  Let \begin{equation}\label{cantor_X} \mathcal{X} = \{ [a,b] : a,b \in \R\}\end{equation} be the collection of closed intervals, $m=2$, and let $M = [0,1]^2 \subseteq \R$. For each $\mathbf{k}=(k^{(1)},k^{(2)}) \in M$, we consider the following two maps, 
	\begin{eqnarray*}
		f^{(1)}_\mathbf{k} : & \mathcal{X} &\to \mathcal{X} \\
		&[a,b]   &\mapsto  [a, k^{(1)}(b-a)+a]
	\end{eqnarray*}
	\begin{eqnarray*}
		f^{(2)}_\mathbf{k} : &\mathcal{X} &\to \mathcal{X} \\
		& [a,b]   &\mapsto [k^{(2)}(a-b)+b, b].
	\end{eqnarray*}

	Note that both $f^{(1)}_\k$ and $f^{(2)}_{\k}$ are compression maps for any $\k \in M$.  Thus, this defines a marking
	\begin{eqnarray*}
		\F : &M &\to {\color{blue}\mathcal{C}(\mathcal{X})^2} \\
		&\k   &\mapsto  f_{\k}= (f^{(1)}_{\k}, f^{(2)}_{\k}).
	\end{eqnarray*}
 }
	Since
 \[\diam \left(f^{(i)}_p([a,b])\right)=p_i \cdot \diam([a,b]) \quad \text{ for } i = 1, 2,\]
 we get that $L \left(f^{(i)}_p \right)=p_i=U\left(f^{(i)}_p\right)$, and hence \begin{equation}\label{cantor_Lk} \mathbf{L}_p=p=\mathbf{U}_p.\end{equation}

Let $E_0 = [0,1] \in \mathcal{X}$ be fixed. 
For any sequence $\{\k_{\ell}\}_{\ell=0}^{\infty} \in M$, towards using Definition~\ref{def: revised} for our $\mathcal{F}$-limit sets, we define the following:
	\begin{eqnarray*}
		&E^{(0)} &= E_0\\
		&E^{(1)} &= f^{(1)}_{\mathbf{k}_0}  \left( E_0 \right) \cup  f^{(2)}_{\mathbf{k}_0}  \left( E_0 \right) =: E_{1} \cup E_{2}\\
		&E^{(2)} &= f^{(1)}_{\mathbf{k}_1}  \left( E_{1} \right) \cup  f^{(2)}_{\mathbf{k}_1}  \left( E_{1} \right) \cup  f^{(1)}_{\mathbf{k}_2}  \left( E_{2} \right) \cup  f^{(2)}_{\mathbf{k}_2}  \left( E_{2} \right) \\
		& &:=\qquad E_{3} \qquad  \cup \quad  E_{4} \quad \cup \quad  E_{5} \quad  \cup \quad  E_{6}\\
		&  \vdots \\
		&E^{(n)}&= \bigcup_{i=2^{n-1}-1}^{2^n-2} \left(f_{\k_i}^{(1)}(E_i) \cup f_{\k_i}^{(2)}(E_i) \right):=\bigcup_{i=2^{n-1}-1}^{2^n-2} \left(E_{2i+1} \cup E_{2i+2} \right)=\bigcup_{ \ell= 2^n-1}^{2(2^n-1)} E_{\ell}.
	\end{eqnarray*}
Observe that the process of constructing the sequence  $\{E^{(n)}\}_{n=0}^{\infty}$ here is independent of the values of $\{\k_{\ell}\}_{\ell=0}^{\infty}$.  For the constant sequence $\mathbf{k}_{\ell}=(\frac{1}{3}, \frac{1}{3})$ for all $\ell$, $E^{(n)} $ is the $n^{th}$-generation of the Cantor set $\C$ and the $\mathcal{F}$-limit set $\displaystyle F = \cap_{n} E^{(n)} =\C $. 

To allow for more general outcomes, we can update the linear functions $f^{(1)}_\k$ and $f^{(2)}_\k$ simply by changing the value of $\k$ at each stage of the construction, which does not change the computational complexity of the process.  
Using this idea, we now construct some examples of Cantor-like sets by choosing suitable sequences $\{\k_{\ell}\}_{\ell=0}^{\infty}$.
	
	\begin{figure}[h] 
		\centering
		\includegraphics[width=4in]{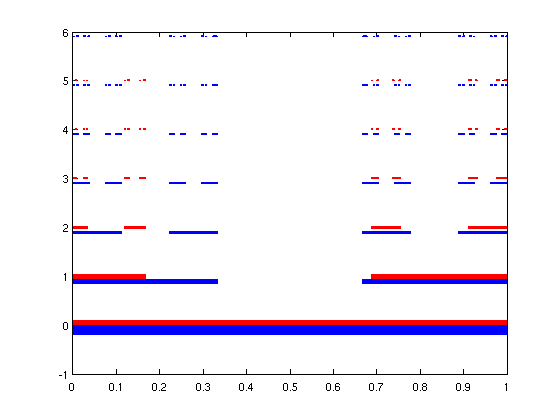} 
		\caption{Comparison of classical Cantor set (blue) and new Cantor-like set (red) }
		\label{comparecant}
	\end{figure}

	\begin{example}\label{example: 1}
		Let $\mathbf{k}_{\ell} = \left(\frac{\ell+1}{4\ell+6}, \frac{2\ell +5}{8\ell+16}\right)$ for $\ell \geq 0$, and let $F$ be the $\F$-limit set generated by the sequence $\{\k_{\ell} \}_{\ell =0}^{\infty}$ with initial set $E_0$. 
  In Figure \ref{comparecant} we plot the usual Cantor set $\C$ (in blue) below the set $F$ (in red) to illustrate the comparison. We can see that the set $F$ has the same basic shape as the Cantor set $\C$, but is no longer strictly self-similar.
		In order to compute the Hausdorff dimension of the new Cantor-like set $F$, we apply Corollary \ref{cor:limit_sup_inf}.  Note that by equation (\ref{cantor_Lk}),
		\[\lim_{\ell \rightarrow \infty}||\mathbf{L}_{\mathbf{k}_\ell}||_s=\lim_{\ell \rightarrow \infty}||\mathbf{k}_\ell ||_s=\frac{2^{\frac{1}{s}}}{4}.\]
		So,  
		\[
		\underline{s_*} =\sup_s\{\liminf_{\ell \rightarrow \infty}||\mathbf{L}_{\mathbf{k}_\ell}||_s>1\} = \sup_s \left\{\frac{2^{\frac{1}{s}}}{4}>1\right\}=\frac{1}{2}.
		\]
		Similarly, we also have $\overbar{s^*}=\frac{1}{2}$. By the following Proposition \ref{UCC_Cantor}, since \[\sup \left \{k_{\ell}^{(1)}+k_{\ell}^{(2)}: \ell=0,1,2, \cdots \right \}=\frac{1}{2}<1,\] $F$ satisfies the uniform covering condition (\ref{condtion_lower}).
		 By Corollary \ref{cor:limit_sup_inf}, $\dim_H(F)=\frac{1}{2}$. 
	\end{example}

 \begin{proposition}
	\label{UCC_Cantor}
	Let $\{\k_{\ell}\}_{\ell=0}^{\infty}$ be a sequence in $M$ with
 \begin{equation}
 \label{eqn: k_upperbound}
     \sup \left \{k_{\ell}^{(1)}+k_{\ell}^{(2)}: \ell=0,1,2, \cdots \right \}<1,
 \end{equation}
	and $F$ be the $\F$-limit set generated by the sequence $\{\k_{\ell} \}_{\ell =0}^{\infty}$ with initial set $E_0$. 
	Then $F$ satisfies the uniform covering condition (\ref{condtion_lower}).
\end{proposition}	
\begin{proof}
Define $\{J_\sigma : \sigma \in D\}$ as in (\ref{F-limit-set}), and let 
\begin{equation} \label{cantor_lambda} 
 \gamma :=\inf_{\ell} \left \{1-k_{\ell}^{(1)}-k_{\ell}^{(2)} \right \}.
 \end{equation}
 By  (\ref{eqn: k_upperbound}), $\gamma\in (0,1]$.
	We will show that $F$ satisfies the uniform covering condition by showing that for any closed interval $B$ in $\mathbb{R}$ with $B\cap F \neq \emptyset$, there exists a $\sigma^*\in D$ such that
$B\cap F\subseteq J_{\sigma^*}$
and $diam(B)\ge \gamma \cdot diam(J_{\sigma^*})$.

Indeed, consider the set
	\[\mathcal{L}:=\{\ell(\sigma): B\cap F \subseteq J_{\sigma}, \sigma\in D\},\]
	where $\ell(\sigma)$ is given in (\ref{ell_sigma}).
	Note that $\mathcal{L}$ is nonempty because $B\cap F \subseteq J_{\emptyset}$ implies that $\ell(\emptyset)\in \mathcal{L}$.
	
 Case 1: If $\mathcal{L}$ is an infinite set, then since $diam(J_{\sigma})\rightarrow 0$ as $\ell(\sigma)\rightarrow \infty$, there exists $\sigma^*\in D$ such that $\ell(\sigma^*)\in \mathcal{L}$ and $diam(B)\ge diam(J_{\sigma^*})\ge \gamma \cdot diam(J_{\sigma^*})$.
	
 Case 2: If $\mathcal{L}$ is finite, let $\ell(\sigma^*)$ be the maximum number in $\mathcal{L}$ for some $\sigma^*\in D$.
	Then,  $\ell(\sigma^*)\in \mathcal{L}$ but $\ell(\sigma^* * j)\notin \mathcal{L}$ for each $j=1,2$. This implies that  $B\cap J_{\sigma^* *j }\neq \emptyset$ for both $j=1, 2$ because $ J_{\sigma^* }= J_{\sigma^* *1 }\cup  J_{\sigma^* *2} $. Since $B$ is an interval, the gap $J_{\sigma^*}\setminus (J_{\sigma^* * 1}\cup J_{\sigma^* * 2})$ between $J_{\sigma^* * 1}$ and $J_{\sigma^* * 2}$ is contained in $B$, which yields that
	\begin{eqnarray*}
		diam(B)&\ge& diam \left (J_{\sigma^*}\setminus (J_{\sigma^* * 1}\cup J_{\sigma^* * 2})\right )\\
		&=& diam(J_{\sigma^*})-diam (J_{\sigma^* * 1})-diam (J_{\sigma^* * 2})\\
		&\ge&  diam(J_{\sigma^*}) \left (1-k_{\ell(\sigma^*)}^{(1)}-k_{\ell(\sigma^*)}^{(2)} \right )\ge \gamma \cdot diam(J_{\sigma^*}).
	\end{eqnarray*}
	As a result, in both cases, the uniform covering condition (\ref{condtion_lower}) holds.
\end{proof}
	
	In the next example, we will construct a random Cantor-like set as follows.
	
	\begin{example} \label{example: 2}
 	\begin{figure}[h] 
			\centering
			\includegraphics[width=3.1in]{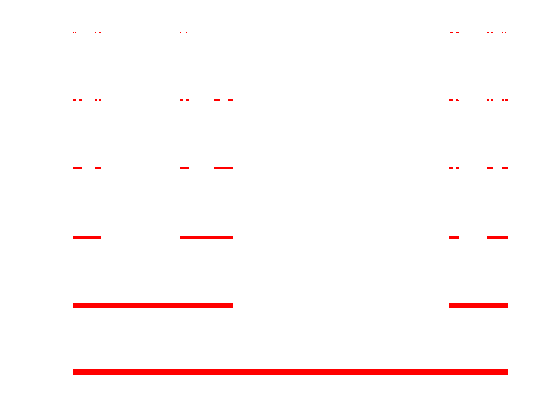} 
			\caption{A randomly generated Cantor-like set}
			\label{cantorrand}
		\end{figure}
		For each $\ell\ge 0$, we take $\mathbf{k}_{\ell}=\left(q_\ell, \frac{1}{2}-q_\ell \right)$ where $q_\ell$ is a random number between $\frac{1}{8}$ and $\frac{3}{8}$. Let $F$ be the corresponding $\F$-limit generated by the sequence $\{\k_{\ell} \}_{\ell =0}^{\infty}$ with initial set $E_0$. We plot the first few generations in Figure \ref{cantorrand}.  In this example, the total length of the $n^{th}$ generation $E^{(n)}$ is chosen to be $(\frac{1}{2})^n$, while the scaling factors of the left subintervals at each stage are randomly chosen.
		
		We now estimate the dimension of $F$.  By (\ref{cantor_Lk}),
		\[\left(\frac{1}{8},\frac{1}{8}\right)\le \mathbf{L}_{\mathbf{k}_\ell} =\mathbf{k}_\ell=\mathbf{U}_{\mathbf{k}_\ell} \le \left(\frac{3}{8},\frac{3}{8}\right).\]
		By Corollary \ref{cor:uniform_bounds}, 
		\[\frac{\log(2)}{-\log(1/8)}\le \dim_H(F) \le \frac{\log(2)}{-\log(3/8)}.\]
		That is,
		\[\frac{1}{3}\le \dim_H(F) \le \frac{\log(2)}{\log(8/3)}\approx 0.7067.\]
		Note that due to Proposition \ref{UCC_Cantor}, $F$ satisfies the uniform covering condition (\ref{condtion_lower}) since $q_{\ell} + (\frac{1}{2}-q_{\ell}) = \frac{1}{2} <1$ for each $\ell \ge 0$.
	\end{example}

	\begin{example} \label{example: 3}
 		\begin{figure}[h] 
			\centering
			\includegraphics[width=3.1in]{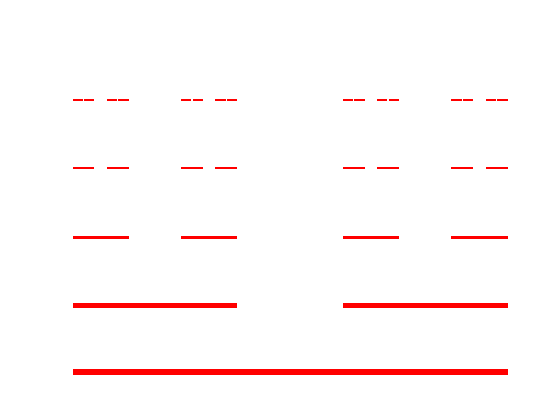} 
			\caption{Fractal of measure $\frac{1}{3}$ created by using $\sum_{n=0}^{\infty} \frac{1}{n!} = e$  }
			\label{example3cant}
		\end{figure}
		In this example, we create a sequence $\{\k_{\ell}\}_{\ell=0}^{\infty}$ that results in a limit set with a given measure, e.g. 1/3.   Of course, the classic example of such a limiting set is the fat Cantor set.  For a different approach, let $\sum_{n=0}^{\infty} a_n$ be any convergent series of positive terms with limit $L$.  We consider a sequence $\{\k_{\ell}\}_{\ell=0}^{\infty}$ defined in the following way.  
		
		Let $n\geq1$ be the generation of the construction and for each $\ell$ with $2^{n-1}-1 \leq \ell \leq 2^{n}-2$, define $\mathbf{k}_\ell=(b_n, b_n)$ where $$\displaystyle b_1:= \dfrac{\frac{3}{2}L-a_0}{2\left(\frac{3}{2}L\right)} \quad \text{ and } \quad b_n := \dfrac{\frac{3}{2}L-\sum_{i=0}^{n-1} a_i}{2\left(\frac{3}{2}L-\sum_{i=0}^{n-2}a_i\right)} \text{ for } n\ge 2.$$ 
		With this sequence $\{\k_{\ell}\}_{\ell =0}^{\infty}$, one can find that the length of each interval in the $n^{th}$ generation is $$b_1 b_2 \cdots b_n=\dfrac{\frac{3}{2}L-\sum_{i=0}^{n-1} a_i}{2^n \cdot \frac{3}{2}L}.$$  Thus, the total length of the $n^{th}$ generation is \[ \dfrac{\frac{3}{2}L-\sum_{i=0}^{n-1} a_i}{\frac{3}{2}L}=1-\frac{2}{3L}\sum_{i=0}^{n-1} a_i \] which converges to 1/3 as desired.  As an example, we take the convergent series $\displaystyle \sum_{n=0}^{\infty} \dfrac{1}{n!} = e$ and use it to create the $\F$-limit set $F$ with measure 1/3.  The first few generations are shown in Figure \ref{example3cant}.
	\end{example}
\end{subsection}

\begin{subsection}{Sierpi\'nski Triangle}
\label{sierpinski examples}
To construct our Sierpi\'nski-like triangles, in this subsection we choose $D$ to be the $3$-ary tree, $M := [0, 1]^6$, and 
\begin{equation}\label{xdefinesierp}
\mathcal{X} := \{(A, B, C) : A, B, C \in \mathbb{R}^2\}
\end{equation}
consists of all triangles $\Delta ABC$ in $\R^2$. 
%
\begin{figure}[h]
		\centering
		\includegraphics[width=3.2in]{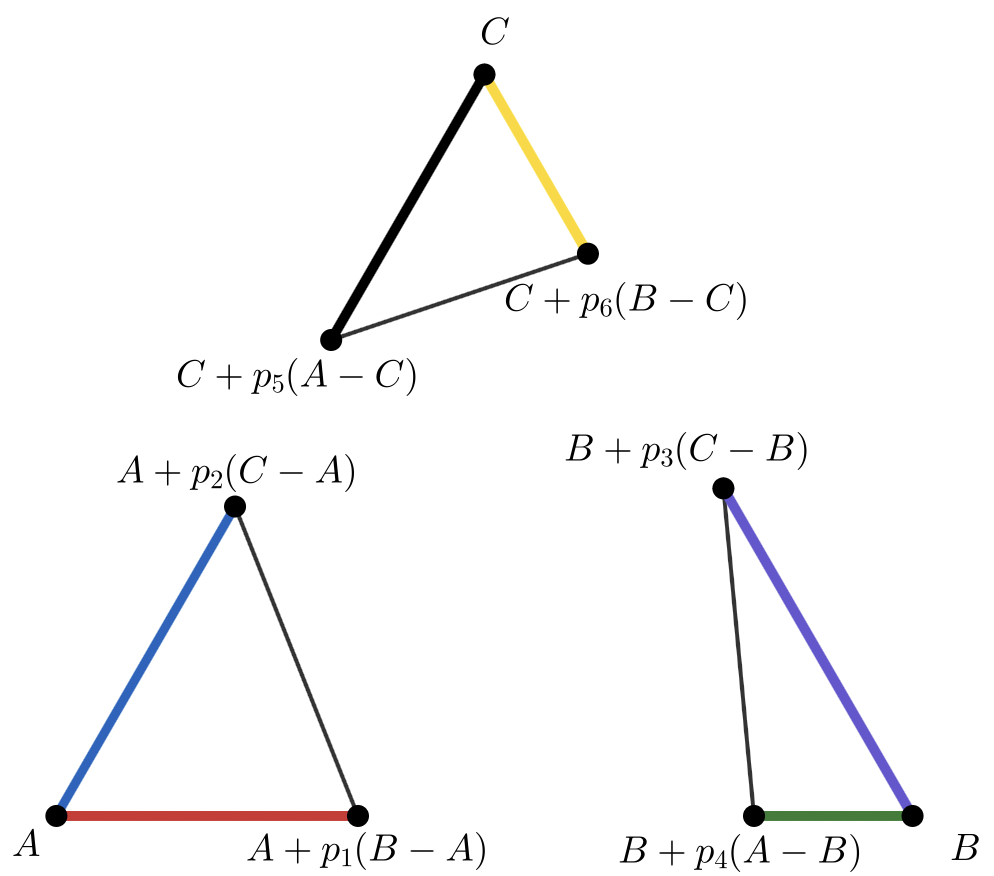}
		\caption{Geometric illustration of $p = (p_1, p_2, p_3, p_4, p_5, p_6) \in M$}
        \label{sierp_lengths}
	\end{figure}
 
Define the marking $\mathcal{F}$ of $M$ into $\mathcal{C}(\mathcal{X})^3$ by 
\[
\mathcal{F}(p) := (f_p^{(1)}, f_p^{(2)}, f_p^{(3)}) \quad \text{ for all } \quad p = (p_1, p_2, p_3, p_4, p_5, p_6) \in [0, 1]^6
\]
where
	\begin{eqnarray*}
		f^{(1)}_{p} (A,B,C) &:= &(A, A+p_1(B-A), A+p_2(C-A))\\
		f^{(2)}_{p} (A,B,C) &:= &(B+p_4(A-B), B, B+p_3(C-B))\\
		f^{(3)}_{p} (A,B,C) &:= &(C+p_5(A-C), C+p_6(B-C), C)
	\end{eqnarray*}
are compression maps from $\mathcal{X}$ to $\mathcal{X}$ as illustrated in Figure~\ref{sierp_lengths}.
	
	Of course, to prevent overlaps we can require that
 \begin{equation}\label{overlap constraints}
 p_1 + p_4 \leq 1, \quad p_2 + p_5 \leq 1, \quad p_3 + p_6 \leq 1.
 \end{equation}
 When each of the inequalities are strict, the images of $f^{(i)}_{p}$ are three disconnected triangles, as illustrated in Figure \ref{disconnect}. When all equalities hold, the images are connected, as illustrated in Figure \ref{connect}. 
	
	\begin{figure}[h]
		\centering
		\begin{subfigure}[b]{0.3\textwidth}
			\includegraphics[width=2in]{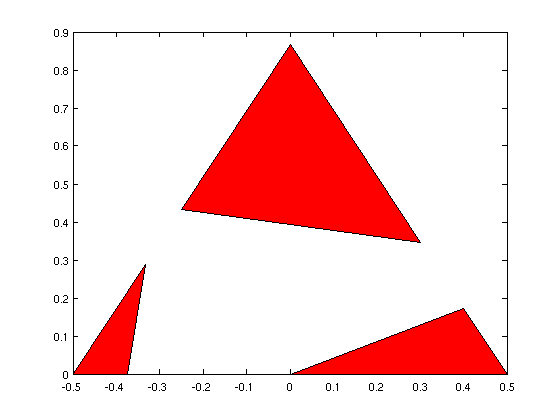}
			\caption{}
			\label{disconnect}
		\end{subfigure}
		~ 
		\qquad \qquad
		\begin{subfigure}[b]{0.3\textwidth}
			\includegraphics[width=2in]{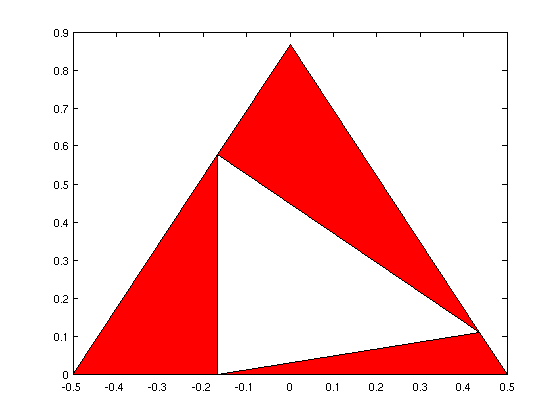}
			\caption{}
			\label{connect}
		\end{subfigure}
		\caption{First generation of disconnected and connected triangles}\label{gen1}
	\end{figure}
	
	In the case of the connected sets, the values of $p = (p_1, p_2, p_3, p_4, p_5, p_6) $ are determined by $p_1, p_2, p_3$ since $p_4 = 1-p_1$, $p_5=1-p_2$, $p_6=1-p_3$.  
 In this case, we may also view $p= (p_1, p_2, p_3)$ as a vector in $[0,1]^3 \subseteq \R^3$.

	To create the normal Sierpi\'nski triangle, we choose
	\begin{equation}
	E_0 = \begin{pmatrix}
	-1/2 & 1/2 & 0\\
	0&0&\sqrt{3}/2
	\end{pmatrix},
	\label{e0sierp}
	\end{equation}
	the equilateral triangle of unit side length, and $\{\k_{\ell}\}_{\ell = 0}^\infty$ to be the constant sequence $\k_{\ell} =\k := (1/2, 1/2, 1/2, 1/2, 1/2, 1/2)$ in $M$ so that each iteration maps a triangle to three triangles of half the side length with the desired translation.  In this case the $\F$-limit set generated by the constant sequence $\{\k_{\ell}=\k \}_{\ell =0}^{\infty}$ with initial set $E_0$  corresponds to the standard Sierpi\'nski Triangle as seen in Figure \ref{Sierpinski triangle}.  
	\begin{figure}[h]
	\centering
	\includegraphics[scale=.18]{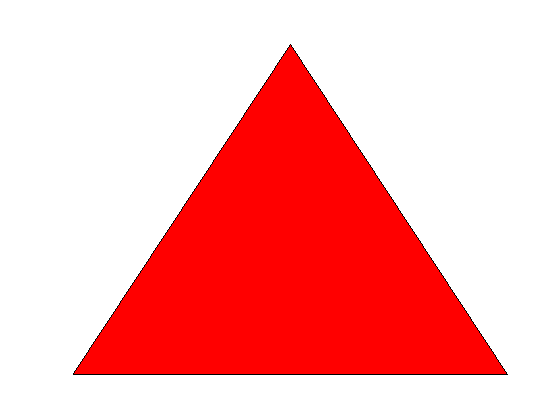}
	\includegraphics[scale=.18]{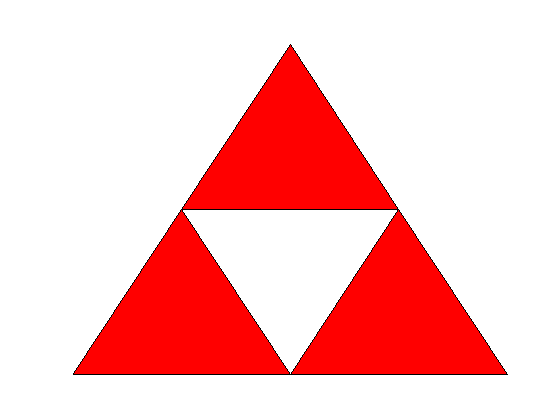}		\includegraphics[scale=.18]{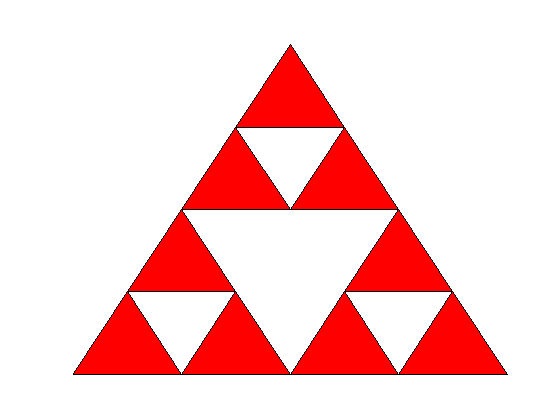} $\cdots$
	\includegraphics[scale=.18]{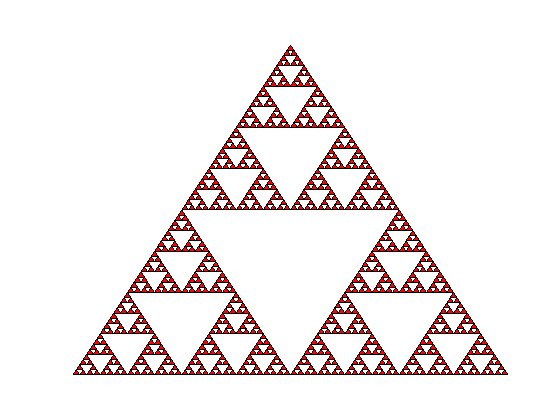}
	\caption{Constructing the Sierpi\'nski triangle}
	\label{Sierpinski triangle}
\end{figure}

To generate Sierpi\'nski-like fractals, we now adjust the values of the marking parameters $\{\k_{\ell}\}_{\ell=0}^{\infty}$.  For each $ p = (p_1, p_2, p_3, p_4, p_5, p_6) \in M=[0, 1]^6$ and $1\le i \le 3$, 
	\begin{eqnarray*}
		U\left(f_{p}^{(i)}\right)=\sup_{(A,B,C)\in \mathcal{X}}\frac{diam \left(f_{p}^{(i)}(A,B,C)\right)}{diam\left((A,B,C)\right)} =\max\left \{p_{_{2i-1}}, p_{_{2i}} \right\},
	\end{eqnarray*}
	and
	\begin{eqnarray*}
		L\left(f_{p}^{(i)} \right)=\inf_{(A,B,C)\in \mathcal{X}}\frac{diam \left(f_{p}^{(i)}(A,B,C)\right)}{diam\left((A,B,C)\right)} =\min\left \{p_{_{2i-1}}, p_{_{2i}} \right\}.
	\end{eqnarray*}

	When $p\in M$ is bounded, i.e. if $0<\lambda \le p_j \le \Lambda<1$ for all $j=1,\cdots ,6$, then by also considering the constraints~(\ref{overlap constraints}),
	\[\mathbf{U}_p \le \mathbf{r}:=(r,\cdots,r) \text{ and } \mathbf{L}_p\ge \mathbf{s}:=(s,\cdots,s),\]
 where $r=\min\{1-\lambda, \Lambda\}$ and $s=\max\{1-\Lambda, \lambda\}$. 
 We usually set $\lambda+\Lambda=1$ so that $r=\Lambda$ and $s=\lambda$.

	Following our general process, we construct some random but connected Sierpi\'nski-like sets by introducing randomness into the choice of the sequence $ \{\k_{\ell}\}_{\ell=0}^{\infty}$ in $[0, 1]^3$.
	\begin{example}\label{example: 4}
		Let $\{\k_{\ell}\}_{\ell=0}^{\infty} =\left \{\left (k_\ell^{(1)}, k_\ell^{(2)}, k_\ell^{(3)} \right )\right \}_{\ell=0}^{\infty}$ be a sequence in $[0,1]^3$ with each $k_{\ell}^{(i)}$ a random number between given numbers $\lambda$ and $\Lambda$ for each $i=1,2,3$. Let $F$ be the $\F$-limit set generated by the sequence $\{\k_{\ell} \}_{\ell =0}^{\infty}$ with initial set $E_0$. Then the $6^{th}$ generation of the construction results in images like Figure \ref{rand1}. Here, in Figure \ref{rand1a}, $\lambda=\frac{1}{4}$ and $\Lambda=\frac{3}{4}$;  while in Figure \ref{rand1b}, $\lambda=0.45$ and $\Lambda=0.55$. Note that the sets are no longer self-similar.
	\end{example}
	\begin{figure}[h]
		\centering
		\begin{subfigure}[b]{0.42\textwidth}
			\includegraphics[width=2.5in]{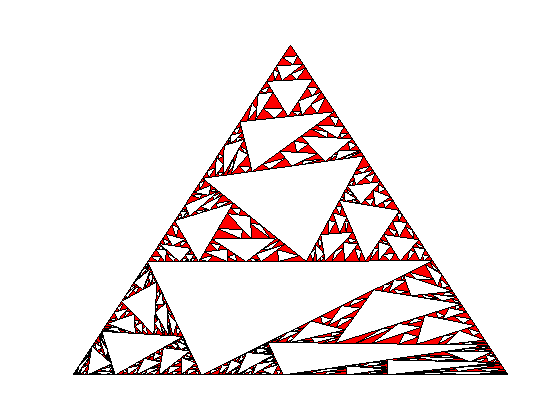}
			\caption{Each $ k_{\ell}^{(i)}$ is random in $[\frac{1}{4}, \frac{3}{4}]$.}
			\label{rand1a}
		\end{subfigure}
		\qquad \qquad
		\begin{subfigure}[b]{0.42\textwidth}
			\includegraphics[width=2.5in]{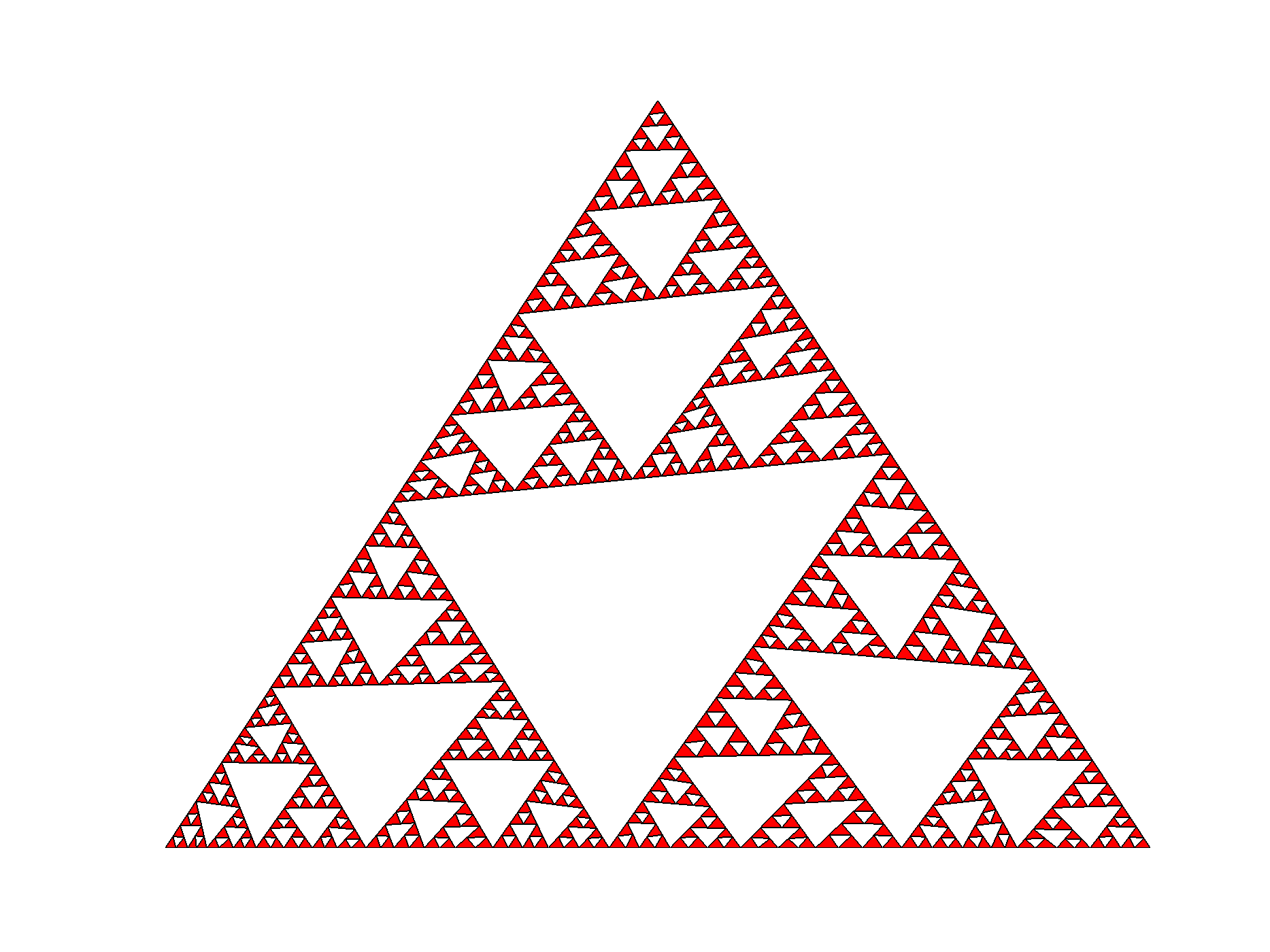}
			\caption{Each $ k_{\ell}^{(i)}$ is random in $[0.45, 0.55]$.}
			\label{rand1b}
		\end{subfigure}
		\caption{Generation 6 of a random and connected Sierpi\'nski triangle}\label{rand1}
	\end{figure}

	In Figure \ref{rand1b}, we pick $\lambda=0.45$ and $\Lambda=0.55$. By Corollary \ref{cor:uniform_bounds},
	\[\frac{\log(m)}{-\log(s)}\le \dim_H(F) \le \frac{\log(m)}{-\log(r)},\]
	where $m=3$, $r=0.55$ and $s=0.45$. That is,
	$1.3758 \le \dim_H(F)\le 1.8377.$
 Here, one may use Theorem~\ref{thm: sufficient_uniform_covering} to show $F$ satisfying the uniform covering condition (\ref{condtion_lower}).

	\begin{example}\label{example: 5}
 \begin{figure}[h] 
		\centering
		\includegraphics[width=3in]{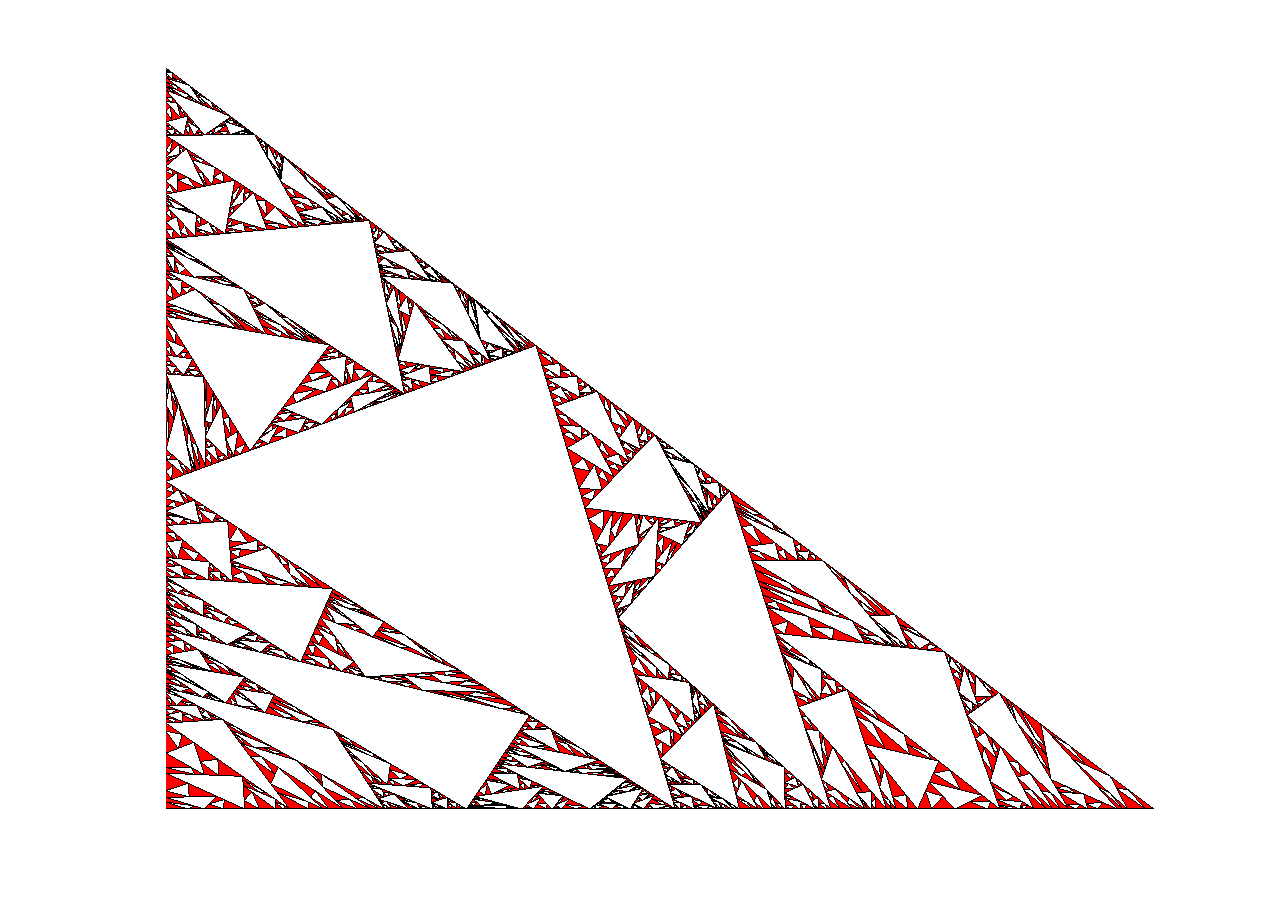} 
		\caption{Generation 7 of a random Sierpi\'nski triangle}
		\label{rand2}
	\end{figure}	
		As in Example \ref{example: 4}, but replacing $E_0$ with  $\tilde{E_0} = \begin{pmatrix} 0 & 1 & 0\\0 & 0&1\end{pmatrix}$, the $7^{th}$ generation of the construction results in an image like Figure \ref{rand2}, when $\lambda=\frac{1}{4}$ and $\Lambda=\frac{3}{4}$.
					
	\end{example}

	\begin{example}\label{example: 5new}  
		Let $\{\k_{\ell}\}_{\ell=0}^{\infty} =\left \{\left (k_\ell^{(1)}, k_\ell^{(2)}, k_\ell^{(3)} \right )\right \}_{\ell=0}^{\infty}$ be a sequence in $[0,1]^3$ where
		\begin{eqnarray*}
			k_{\ell}^{(1)} := \frac{1}{2}+\frac{a_\ell}{\sqrt{\ell+1}},
			\quad k_{\ell}^{(2)} := \frac{1}{2}+\frac{b_\ell}{\sqrt{\ell+1}}, 
			\quad k_{\ell}^{(3)} := \frac{1}{2}+\frac{c_\ell}{\ell+1}
   \end{eqnarray*}
		for random numbers $a_\ell, b_\ell, c_\ell \in [-\frac{1}{3}, \frac{1}{3}]$. Let $F$ be the $\F$-limit set generated by the sequence $\{\k_{\ell} \}_{\ell =0}^{\infty}$ with initial set $E_0$. 
		Then the seventh generation of the construction of $F$ results in an image like Figure \ref{rand3}. 
		
		\begin{figure}[h] 
			\centering
			\includegraphics[width=3in]{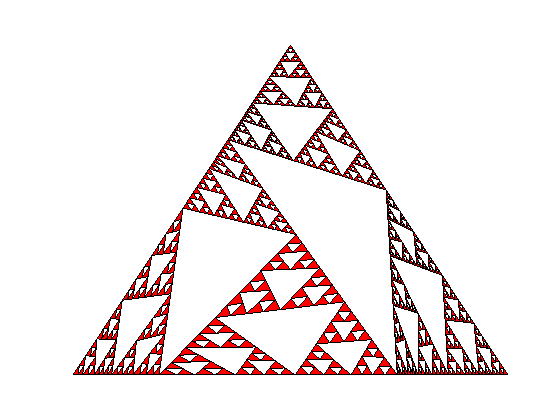} 
			\caption{Generation 6 of a Sierpi\'nski-type triangle with controlled dimension}
			\label{rand3}
		\end{figure}
		
		In this case, we can calculate the exact value of the Hausdorff dimension of $F$.  Indeed, by Corollary \ref{cor:limit_sup_inf},
		\[\lim_{\ell \rightarrow \infty} (||\mathbf{U}_{\mathbf{k}_\ell}||_s)^s=\frac{3}{2^s}=\lim_{\ell \rightarrow \infty} (||\mathbf{L}_{\mathbf{k}_\ell}||_s)^s.\]
		Thus, $\dim_H(F)=\frac{\log(3)}{\log(2)}$, where one may use Theorem~\ref{thm: sufficient_uniform_covering} to show that $F$ satisfies the uniform covering condition (\ref{condtion_lower}).
	\end{example}
	
\end{subsection}

\begin{subsection}{Menger Sponge}
	\label{menger section}
	
	Let \begin{equation} \mathcal{X} = \left \{ (O,A,B,C) : O,A,B,C \in \R^3 \right \} \label{xdefinemeng} \end{equation} representing the collection of all rectangular prisms $(OABC)$ in $\R^3$, $m=20,$ and	
 \begin{equation}\label{3d_M} 
 M = \left \{(p_1, p_2, p_3, p_4, p_5, p_6) \in [0,1]^6 :  p_1 \leq p_2,\ p_3 \leq p_4,\ p_5 \leq p_6 \right \}.
 \end{equation}
	
			\begin{figure}[h]
		\centering
		\includegraphics[width=2.5in]{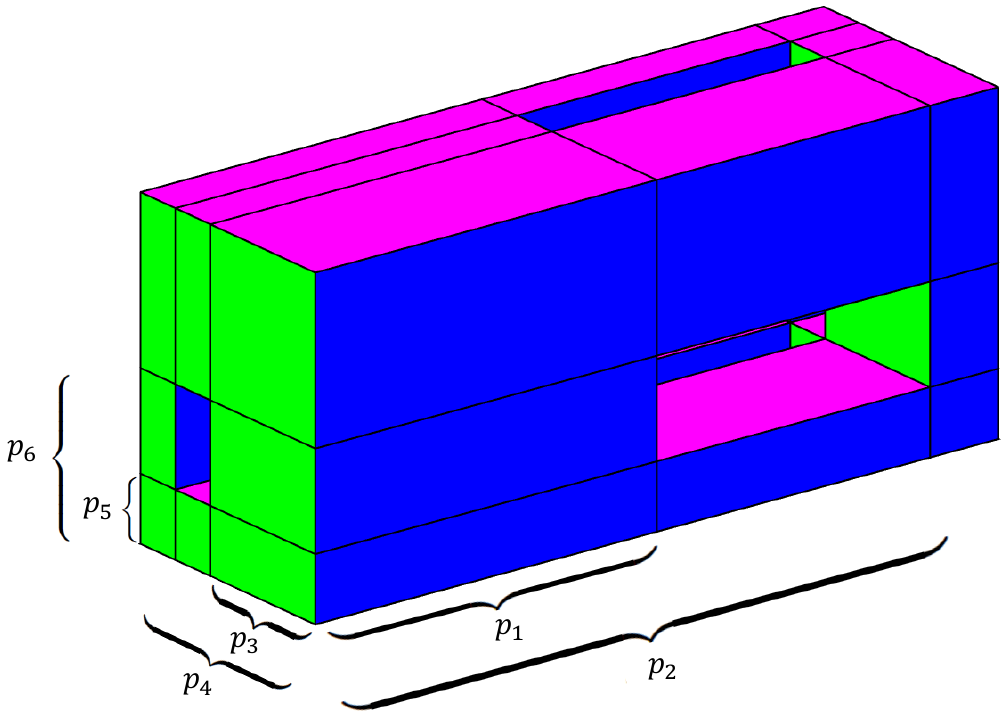}
		\caption{Geometric illustration of $p=(p_1, p_2, p_3, p_4, p_5, p_6) \in M$ used in constructing $f^{(i)}_p$}
		\label{meng_length}
	\end{figure}

	For each $p \in M$ and $i =1, 2, \dots , 20$, we can define affine transformations $f^{(i)}_p : \mathcal{X} \to \mathcal{X}$ as follows.  
 For any $p=(p_1, p_2, p_3, p_4, p_5, p_6) \in M$, define  
	$$ T = \begin{bmatrix} 0& p_1& p_2& 1 \end{bmatrix}, \quad R= \begin{bmatrix} 0& p_3& p_4& 1 \end{bmatrix}, \quad S= \begin{bmatrix} 0& p_5& p_6& 1 \end{bmatrix}.$$
	Let
	\begin{eqnarray*}
		I&=& \{ (a,b,c) : 1 \leq a,b,c\leq3 \text{ with } a,b,c \in \Z, \text{and no two of $a,b,c$ equal to 2} \}.
	\end{eqnarray*}
	For each $(a,b,c) \in I$ and $p \in M$, define 
	\begin{equation*}
	Q_{p}(a,b,c)=\begin{bmatrix} 
	1-( T(a) + R(b) +S(c)) & T(a) & R(b) & S(c)\\
	1-( T(a+1) + R(b) +S(c)) & T(a+1) & R(b) & S(c)\\
	1-( T(a) + R(b+1) +S(c)) & T(a) & R(b+1) & S(c)\\
	1-( T(a) + R(b) +S(c+1)) & T(a) & R(b) & S(c+1)\\
	\end{bmatrix}
	\end{equation*}
 where $T(a)$ denotes the $a^{th}$ entry of the vector $T$, and similarly for the others.
 
	Note that the set $I$ contains 20 elements, so we can express it as $$I = \{ (a_i, b_i, c_i) : 1 \leq i \leq 20\}.$$
	For each $p \in M$ and $1 \leq i \leq 20$, we consider the affine transformation $f^{(i)}_{p} : \mathcal{X} \to \mathcal{X}$ given by
	\begin{equation} 
	f^{(i)}_{p} (O,A,B,C)  = Q_{p} (a_i,b_i,c_i)  \begin{bmatrix} O\\A\\B\\C\end{bmatrix} 
	\label{fdefinemeng_f} 
	\end{equation}
		for every $(O, A, B, C) \in \mathcal{X}$.  
		Note that for $i=1, \dots, 20$ and $p \in M$, $f^{(i)}_{p}$ is a compression.   Thus, we can define a marking $\F : M \to \mathcal{C}(\mathcal{X})^{20}$ by sending $p   \mapsto  f_{p}= (f^{(1)}_{p}, \dots , f^{(20)}_{p})$.
		Using this, for any starting rectangular prism $E_0 = (O,A,B,C) \in \mathcal{X}$, we can generate a sequence of sets that follows a similar construction to the Menger sponge.
	
		In the following examples, we will construct $\mathcal{F}$ limit sets with the unit cube
		\begin{equation}
		E_0 = \begin{bmatrix}
		0 & 1 & 0 & 0\\
		0&0&1&0\\
		0&0&0&1
		\end{bmatrix}
		\label{e0meng}
		\end{equation}
  as the initial set. 
  Note that, for $\k=(1/3, 2/3, 1/3, 2/3, 1/3, 2/3)$,  the $\mathcal{F}$ limit set generated by the constant sequence $\{\k_\ell=\k\}_{\ell=0}^\infty$ with initial set $E_0$ is the classical Menger sponge as shown in Figure \ref{Menger cube}. 
		
\begin{figure}[h]
	\centering
	\includegraphics[scale=.18]{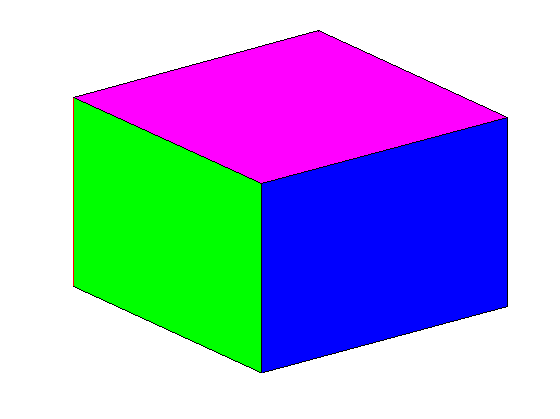}	\includegraphics[scale=.18]{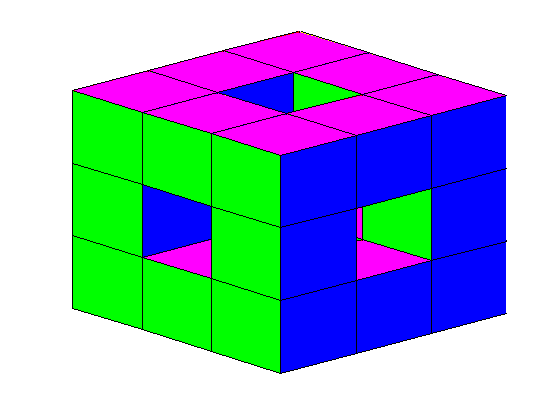}
	\includegraphics[scale=.18]{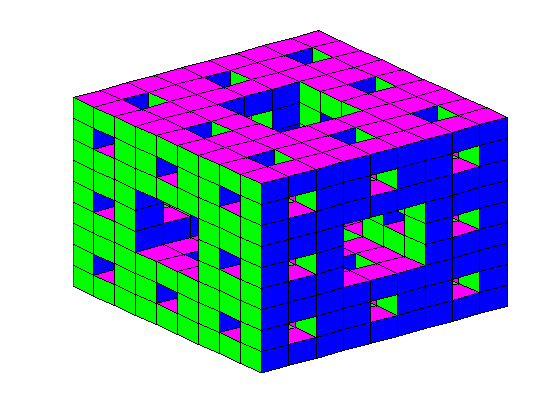}
	\includegraphics[scale=.18]{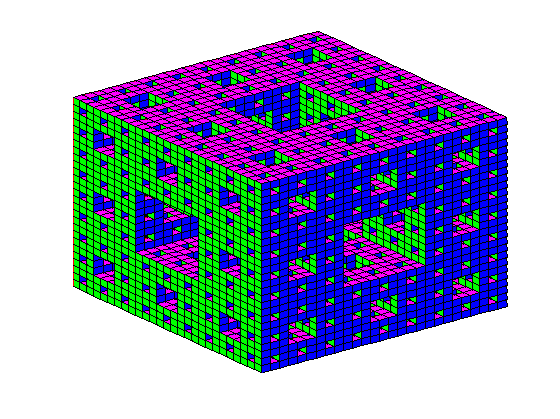} 
	\caption{The first three generations of the Menger sponge}
	\label{Menger cube}
\end{figure}

To estimate the dimension of variations of the Menger sponge, let us make the following calculations. 
For each $p=(p_1, p_2 \cdots, p_6)\in M$ and $1\le i \le 20$, 
	\begin{eqnarray*}
		U\left(f_{p}^{(i)}\right)&=&\sup_{(O,A,B,C)\in \mathcal{X}}\frac{diam \left(f_{p}^{(i)}(O,A,B,C)\right)}{diam\left((O,A,B,C)\right)} \\
		&=&\sup_{(O,A,B,C)\in \mathcal{X}}\frac{diam \left(Q_{p}(a_i, b_i, c_i)[O,A,B,C]'\right)}{diam\left((O,A,B,C)\right)} \\
		&=&\max\{T(a_{i+1})-T(a_i), R(b_{i+1})-R(b_i), S(c_{i+1})-S(c_i)\}.
	\end{eqnarray*}
	Similarly,
	\begin{eqnarray*}
		L\left(f_{p}^{(i)}\right)=\min\{T(a_{i+1})-T(a_i), R(b_{i+1})-R(b_i), S(c_{i+1})-S(c_i)\}.
	\end{eqnarray*}
	When $p_{2j} = 1-p_{2j-1}$ for each $j=1, 2, 3$, it is easy to check that
	\begin{eqnarray*}
		\sum_{i=1}^{20}U(f_{p}^{(i)})^s&=&\sum_{i=1}^{20}\max\{T(a_{i+1})-T(a_i), R(b_{i+1})-R(b_i), S(c_{i+1})-S(c_i)\}^s\\
		&=&8\max\{p_1, p_3, p_5\}^s+4\max\{1-2p_1, p_3, p_5\}^s  \\
		& & +4\max\{p_1, 1-2p_3, p_5\}^s +4\max\{p_1, p_3, 1-2p_5\}^s.
	\end{eqnarray*}
 Note that, when $0<\lambda \le p_1, p_3, p_5 \le \Lambda<1$, it follows that
	\begin{equation}\label{eq: upper bound estimate}
	(||\mathbf{U}_p||_s)^s=\sum_{i=1}^{20}U\left(f_{p}^{(i)}\right)^s \le 8\Lambda^s +12\max\{1-2\lambda, \Lambda\}^s.
	\end{equation}
	Similarly,
	\begin{equation}\label{eq: lower bound estimate}
		(||\mathbf{L}_p||_s)^s \ge 8\lambda^s +12\min\{1-2\Lambda, \lambda\}^s. 
	\end{equation}
 
	\begin{example}\label{example: 9}
		Let $\k_{\ell} =\left(k_\ell^{(1)}, k_\ell^{(2)}, k_\ell^{(3)}, k_\ell^{(4)}, k_\ell^{(5)}, k_\ell^{(6)}\right) \in M$ with each $k_\ell^{(2j-1)}$ a random number between given parameters $\lambda$ and $\Lambda$ and $k_\ell^{(2j)}=1-k_\ell^{(2j-1)}$ for each $j=1,2,3$.  Let $F$ be the $\F$-limit set generated by the sequence $\{\k_{\ell}\}_{\ell =0}^{\infty}$ with initial set $E_0$.
  Then the third iteration of the construction of $F$ results in images like Figure \ref{meng3}.
  Here, in Figure \ref{meng3} (A) the parameters $\lambda=0$ and $\Lambda=\frac{1}{2}$, while in Figure \ref{meng3} (B) the parameters $\lambda=0.32$ and $\Lambda=0.35$.
		
  \begin{figure}[h]
  \centering
  \includegraphics[width=5.1in]{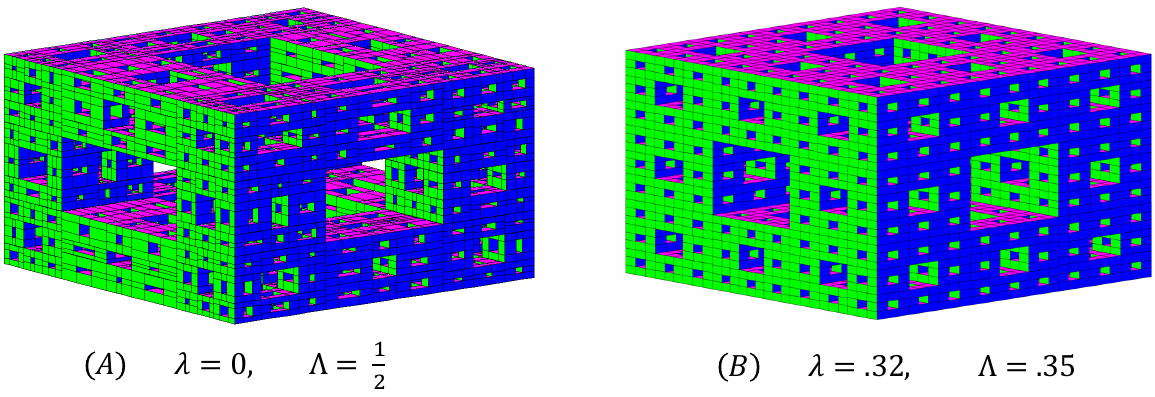}
  \caption{Generation 3 of random Menger sponge}\label{meng3}
  \end{figure}
  
	We now estimate the dimension of the random $\mathcal{F}$ limit set $F$ for $\lambda = .32$ and $\Lambda = .35$ as illustrated in Figure \ref{meng3}(B). 
 By equation~(\ref{eq: upper bound estimate}), for any $s>2.901$,
	\begin{eqnarray*}
		(||\mathbf{U}_{\k_\ell}||_s)^s &\le & 8\Lambda^s+12\max\{1-2\lambda, \Lambda\}^s \le 8*0.35^s+12*0.36^s \\
		&<& 8*0.35^{2.901}+12*0.36^{2.901}\approx 1.000.
	\end{eqnarray*}
	By Theorem \ref{thm:ratio_bounds}, $\dim_H(F)\le 2.901$. 
 Similarly, by equation~(\ref{eq: lower bound estimate}), 
 for any $s \le 2.546$,
	\begin{eqnarray*}
		(||\mathbf{L}_{\k_\ell}||_s)^s
		&\ge & 8\lambda^s+12\min\{1-2\Lambda, \lambda\}^s \\
		&\ge & 8*0.32^s+12*0.3^s \ge 8*0.32^{2.546}+12*0.3^{2.546}\approx 1.000.
	\end{eqnarray*}
	By Theorem \ref{thm:ratio_bounds} again, $\dim_H(F)\ge 2.546$, where one may use Theorem~\ref{thm: sufficient_uniform_covering} to show that $F$ satisfies the uniform covering condition (\ref{condtion_lower}). As a result,
	\[2.546 \le \dim_H(F)\le 2.901.\]
		\end{example}

  In the following example, we are able to calculate the exact Hausdorff dimension of a non self-similar Menger sponge. 
 		\begin{figure}[h] 
			\centering
			\includegraphics[width=3.5in]{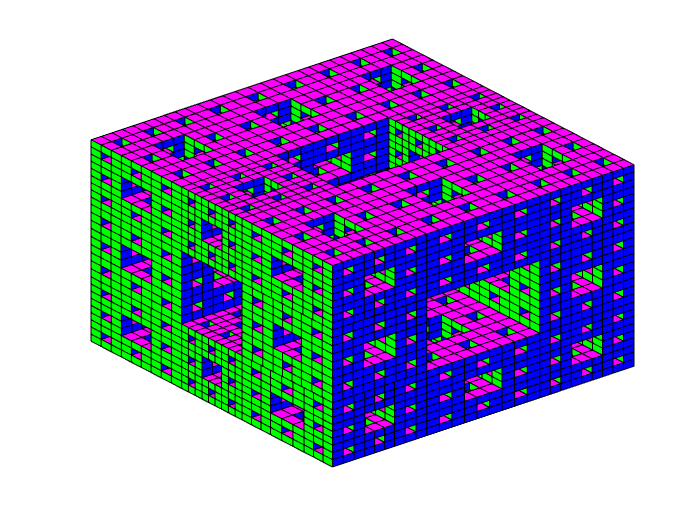} 
			\caption{Generation 3 of a non-self similar Menger sponge with Hausdorff dimension $\frac{\log(20)}{\log(3)}$}
			\label{meng4}
		\end{figure}
	\begin{example}
		\label{3d_exact}
		For each $\ell \ge 0$, let $\mathbf{k}_\ell=\left(k_{\ell}^{(1)}, k_{\ell}^{(2)}, \cdots, k_{\ell}^{(6)}\right)$ where
		\begin{eqnarray*}
			k_{\ell}^{(1)} &=&\frac{1}{3}+\frac{(-1)^\ell}{12(\ell+1)^2},\quad  k_{\ell}^{(2)}=1-k_{\ell}^{(1)},\\
			k_{\ell}^{(3)} &=& \frac{1}{3}-\frac{(-1)^\ell}{6(\ell+1)^2}, \quad k_{\ell}^{(4)}=1-k_{\ell}^{(3)},\\
			k_{\ell}^{(5)} &=& \frac{1}{3}+\frac{(-1)^\ell}{18(\ell+1)^2}, \quad k_{\ell}^{(6)}=1-k_{\ell}^{(5)}.
		\end{eqnarray*}
  Let $F$ be the $\F$-limit set generated by the sequence $\{\k_{\ell} \}_{\ell =0}^{\infty}$ with initial set $E_0$.
	Then the third generation of the construction of $F$ leads to an image like Figure \ref{meng4}.	

	By direct computation, 
	\[\lim_{\ell \rightarrow \infty} (||\mathbf{U}_{\mathbf{k}_\ell}||_s)^s=\frac{20}{3^s}=\lim_{\ell \rightarrow \infty} (||\mathbf{L}_{\mathbf{k}_\ell}||_s)^s.\]
	Thus, by Corollary \ref{cor:limit_sup_inf}, $\dim_H(F)=\frac{\log(20)}{\log(3)} \approx 2.7268$, since $F$ satisfies the uniform covering condition according to Example \ref{3d_ucc}.
 	\end{example}
\end{subsection}

\begin{remark}  
Another way to construct fractals with partial self similarity was introduced by Barnsley, Hutchinson, and Stenflo~\cite{vvar1}. 
These fractals are called $V$-variable fractals, and their dimensions were studied in~\cite{vvar2}. 
Here, we make some comparisons between $\mathcal{F}$-limit sets and $V$-variable fractals.
In essence, the construction of a $V$-variable fractal uses at most $V \in \N$ number of distinct patterns within each generation of the construction.  
This is done through the following process.

Let $(X,d)$ be a metric space, $\Lambda$ an index set, $F^{\lambda} = \{f^{\lambda}_1, f^{\lambda}_2, \dots, f^{\lambda}_m\}$ an IFS for each $\lambda \in \Lambda$, and $P$ a probability distribution on some $\sigma$-algebra of subsets of $\Lambda$.  Then denote $\mathbf{F} = \{(X,d), F^{\lambda}, \lambda \in \Lambda, P\}$ to be a family of IFSs (with at least two functions in each IFS) defined on $(X,d)$.  Assume that the IFSs $F^{\lambda}$ are \textit{uniformly contractive} and \textit{uniformly bounded}, that is, for some $0<r<1$,
\begin{eqnarray}
&&\sup_{\lambda} \max_{m} d \left( f^{\lambda}_m (x), f^{\lambda}_m(y) \right) \leq r d(x,y),\\
&&\sup_{\lambda} \max_m d\left( f^{\lambda}_m(a),a \right)< \infty
\end{eqnarray} for all $x,y \in X$ and some $a \in X$.

A {\it tree code} is a map $\omega : D \to \Lambda$ from a tree $D$ into the index set $\Lambda$.
Given any $\sigma \in D$, one can naturally define another tree code $(\omega \mylee \sigma) (\tau) := \omega(\sigma \ast \tau)$ that ``starts'' at $\sigma$.
A tree code $\omega$ is {\it $V$-variable} if for each positive integer $k$, there are at most $V$ distinct tree codes of the form $\omega\mylee \sigma$ with $\sigma \in D_k$.
For example, consider the Sierpi\'nski triangle.  
We let $F$ be the IFS that maps the triangle to three copies of 1/2 the size, as usual.  
Let $G$ be the IFS that maps the initial triangle to three triangles that are 1/3 the size, with the vertices shared with the initial set being the fixed points of the maps.  
See Figure \ref{vvargenerators} for the image of the initial step of each.  
\begin{figure}[h]
	\centering
	\includegraphics[scale=.2]{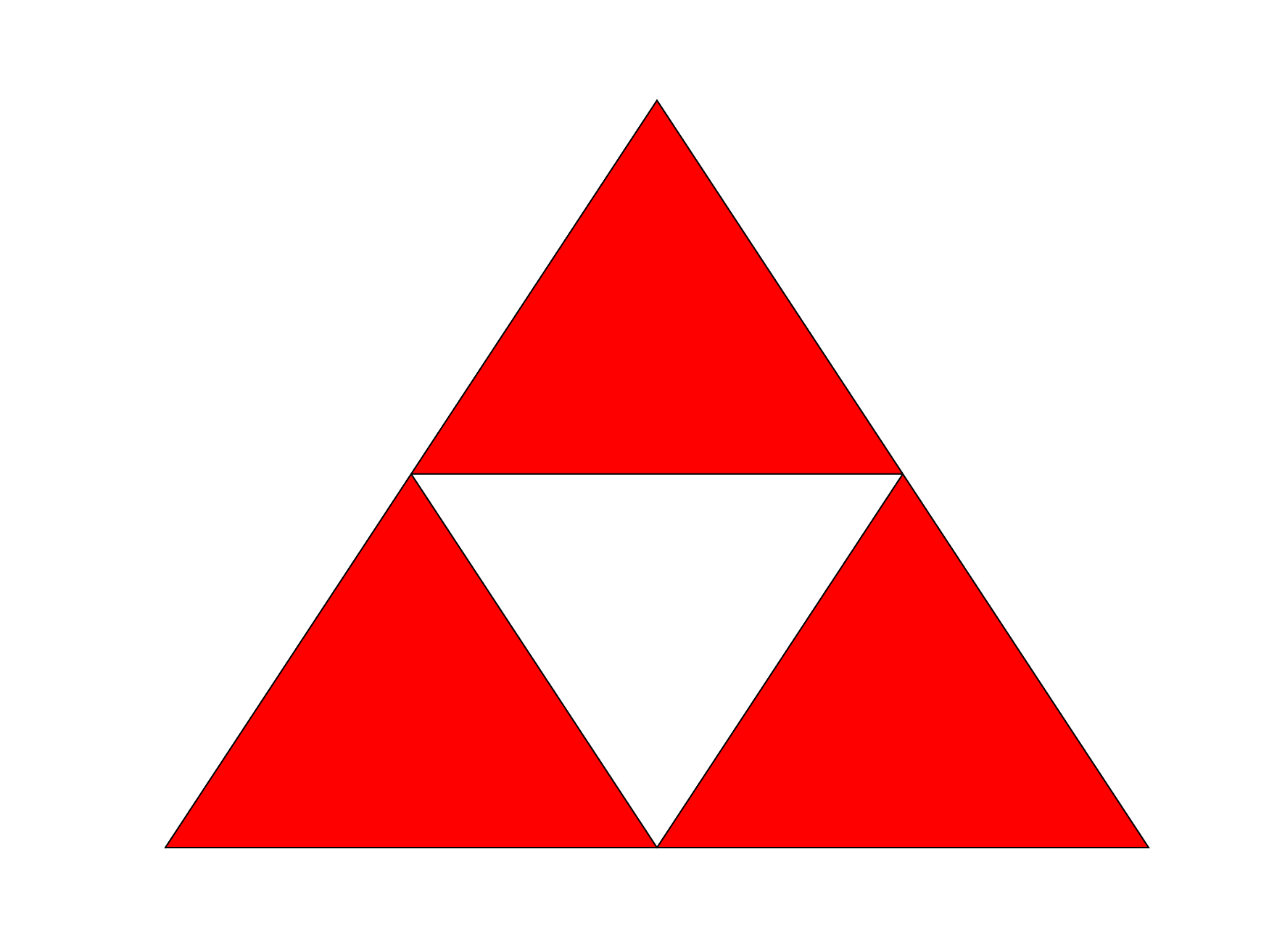}
	\includegraphics[scale=.2]{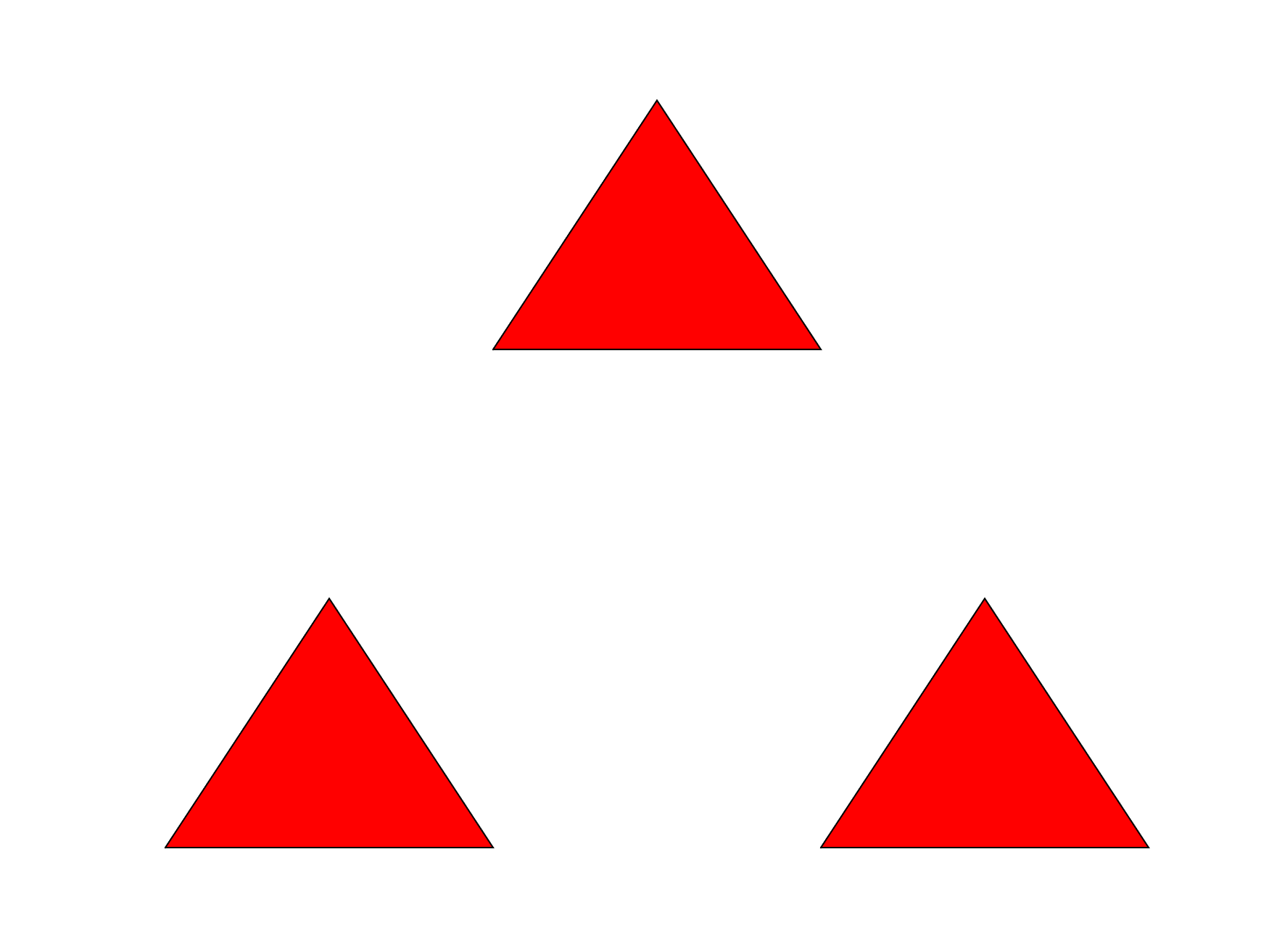}
	\caption{Initial steps of IFSs $F$ and $G$ respectively}
	\label{vvargenerators}
\end{figure}
Thus, $\mathbf{F} = \{ (\R^2, d), \{F,G\}, P=(1/2, 1/2)\}$ is the family $\{F,G\}$ with probability function uniformly choosing 1/2 for each IFS.  
Using these IFSs, three $V$-variable pre-fractals are given in Figure \ref{vvarexample}, being 1-variable, 2-variable, and 3-variable respectively.

Now, we express this $V$-variable fractal in terms of an $\F$-limit set.  
\begin{figure}[h]
	\centering
	\includegraphics[width=3.5in]{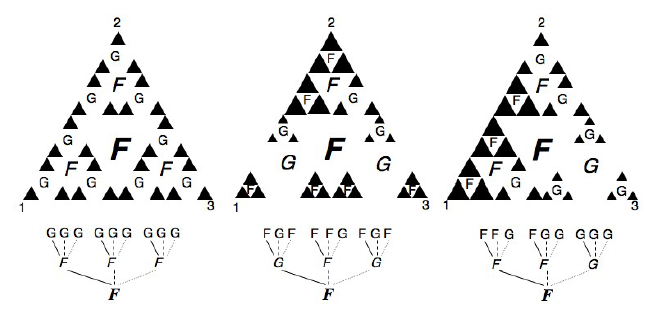}
	\caption{$n=3$ generation prefractals that are 1, 2 and 3-variable respectively.  Images from \cite{vvar2}.}
	\label{vvarexample}
\end{figure}  
Let $\mathcal{X}, M,\mathcal{F}$ and $E_0$ be as in Section \ref{sierpinski examples} and define a map $h: \{F, G\}\rightarrow M$ by $h(F)=(1/2,1/2,1/2,1/2,1/2,1/2)$ and $h(G)=(1/3,1/3,1/3,1/3,1/3,1/3)$. 
For any $V$-variable tree code $\omega : D \to \Lambda = \{F, G\}$, we can define an associated map $\k_\omega : D \to M$ given by $\k_\omega = h\circ \omega$.
Then, the $V$-variable fractal generated by the $V$-variable tree code $\omega$ is the $\mathcal{F}$-limit set generated by $\k_\omega$ with initial set $E_0$.

We now provide some observations of $\mathcal{F}$-limit sets with $V$-variables fractals. First note that the sequence $\k_\omega$ associated with a V-variable tree code $\omega$ contains many repeated terms.
In particular, at any generation $k > V$ there must be at least two nodes $\sigma$ and $\sigma' \in D_k$ such that their corresponding tree-codes are equal, i.e., $\omega\mylee \sigma = \omega\mylee \sigma'$.
This means that $\k_\omega$ is the exact same for a large number of subtrees, whereas in general, as illustrated in Example~\ref{example: 5new}, the maps $\k$ that we use to construct $\mathcal{F}$-limit sets are not. 
In this sense, the $\mathcal{F}$-limit set approach tackles the case of $V$-variable fractals where $V = \infty$.

Here is another observation. To address randomness, we consider maps from $D$ to a parametrization space $M$, while the $V$-variables approach uses a map from $D$ to an index set $\Lambda$. 
In this sense, $M$ plays a similar role as $\Lambda$.
After that, we consider a mapping $\mathcal{F}$ from $M$ to the spaces of all compressions between a collection of sets to itself, while the $V$-variable approach uses IFSs between ambient spaces.

As for the existence of the fractals, in our case, the existence of the limit fractals trivially follows from the definition, while in the case of $V$-variables, one needs to prove the existence of attractors. 
\end{remark}

\section{Sufficient conditions for the Uniform Covering Condition}
\label{UCC}
In previous sections, we have seen that the uniform covering condition (\ref{condtion_lower}) plays a vital role in computing a lower estimate for the Hausdorff dimension of a fractal.  
In this section, Theorem~\ref{thm: sufficient_uniform_covering} provides us with a sufficient condition needed for a fractal to satisfy the uniform covering condition.
This theorem can be used for most examples in \S\ref{sec: examples}. 
For the sake of brevity, we only provide the details for the Menger sponge-like fractal encountered in Example~\ref{3d_exact}. 

Motivated by Proposition \ref{UCC_Cantor}, we introduce two quantities as follows.
The first quantity $\rho_n$ intuitively describes the ``gap'' between $n + 1$ elements of a collection of subsets.
\begin{definition}
Let $n \geq 1$ and $\mathcal{H}$ be a collection of subsets of a metric space $(X,d)$.  Define 
\begin{multline}\label{rho_n}
 \rho_n(\mathcal{H})=\inf\{diam(B) : B \subseteq X \text{ is a closed ball intersecting at least} \\ 
 n+1 \text{ elements in } \mathcal{H}\}.
\end{multline}
\end{definition}
The second quantity $\gamma_n$ describes the minimum ``gap" between any $n+1$ children of a generation relative to the size of their parents.
\begin{definition} Let $\mathcal{J}=\{J_{\sigma} : \sigma \in D\}$ be a collection of compact subsets of a metric space $(X,d)$, and $n \geq 1$.  Define \begin{multline} \label{lambda_n}\gamma_n(\mathcal{J}):=\inf \left \{\frac{\rho_n(\{J_{\sigma*i}: \sigma\in R_k, i=1,2,\cdots, m\})}{\sum_{\sigma \in R_k}diam(J_\sigma)}: \text{ for some $k$} \right.\\ \left. \text{ and } R_k\subseteq D_k \text{ with } 1 \leq |R_k|\le n \right \},\end{multline}  where $|R_k|$ denotes the cardinality of the set $R_k$.
\end{definition}

We now demonstrate how to calculate $\gamma_n(\mathcal{J})$ for the collections $\mathcal{J}$ that we constructed in \S\ref{sec: examples}.
Note that given $\mathcal{X}$, $M$, $D$, and $\mathcal{F}$, for any sequence $\{\k_\ell\}$ in $M$ (or equivalently, any map $\k: D\rightarrow M$), we can construct the collection $\mathcal{J}=\{J_\sigma : \sigma \in D\}$ by (\ref{F-limit-set}).

\begin{example} Let $\mathcal{J}=\{J_{\sigma}: \sigma \in D\}$ be a collection of closed intervals used in the construction of the Cantor-like sets in \S\ref{cantor-like sets}. 
Then 
\begin{eqnarray*}
\gamma_1(\mathcal{J}) &=& \inf \left \{ \frac{\rho_1(\{J_{\sigma*i}: \sigma\in R_k, i=1,2\})}{diam(J_\sigma)}: \text{ for some } R_k\subseteq D_k \text{ with }|R_k|= 1\right \}\\
&=& \inf \left \{ \frac{\rho_1(\{J_{\sigma*1}, J_{\sigma*2}\}) }{diam(J_{\sigma})} : \text{ for } \sigma \in D \right \}\\
&=& \inf \left \{ \frac{diam(J_{\sigma}) - diam(J_{\sigma*1}) - diam(J_{\sigma*2})}{diam(J_{\sigma})}: \sigma \in D \right \}\\
&=& \inf \left \{1-\frac{diam(J_{\sigma * 1})}{diam(J_{\sigma})} - \frac{diam(J_{\sigma * 2})}{diam(J_{\sigma})}: \sigma \in D \right \},
\end{eqnarray*} which agrees with the $\gamma$ in (\ref{cantor_lambda}), see Figure \ref{cantor_rho}.

	\begin{figure}[h]
	\centering
	\includegraphics[scale=.17]{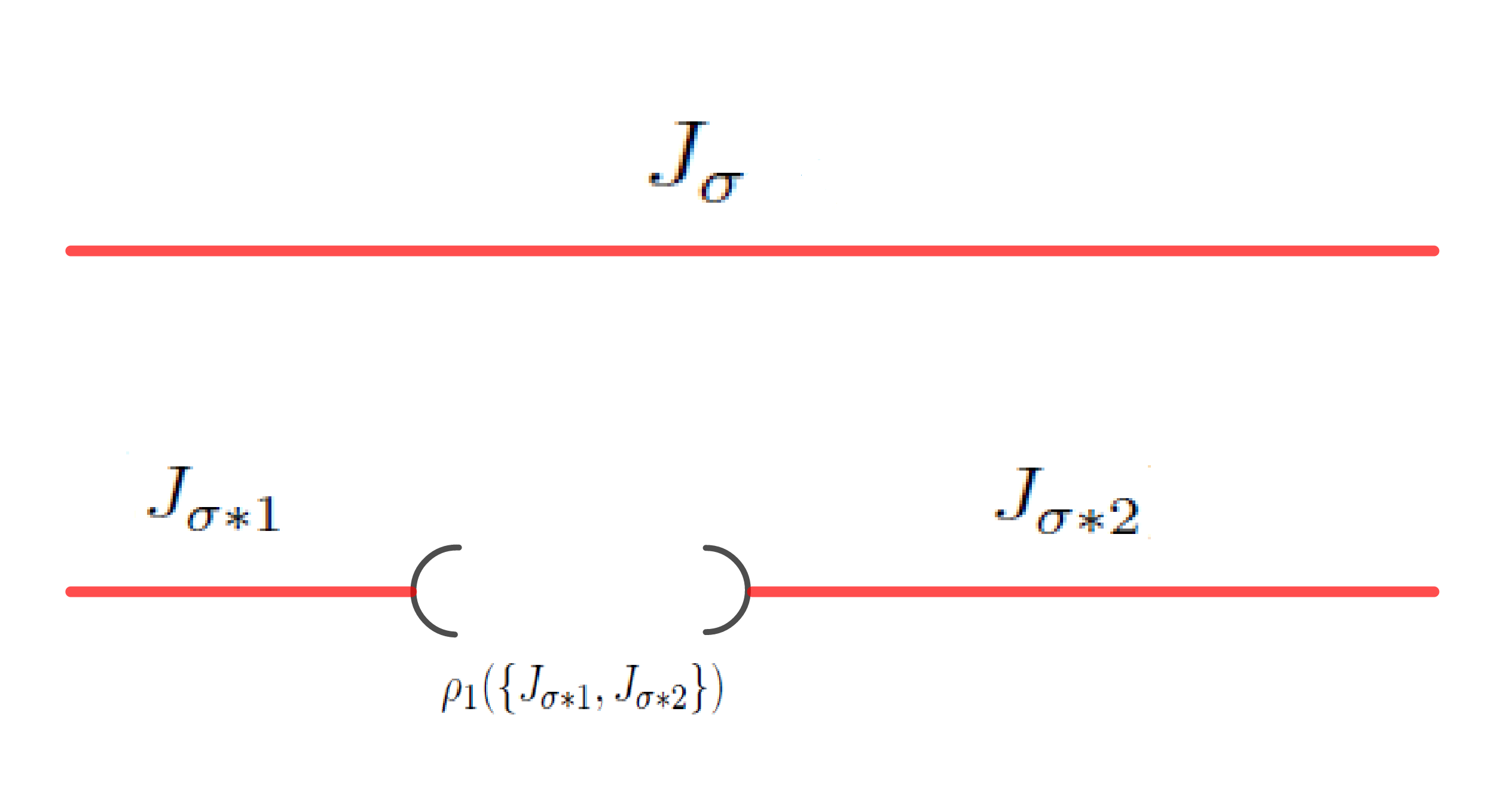}
	\caption{Illustration of $\rho_1(\{J_{\sigma*1}, J_{\sigma*2}\})$}
	\label{cantor_rho}
	\end{figure}
\end{example}

\begin{example}
Let $\mathcal{J}=\{J_{\sigma}: \sigma \in D\}$ be a collection of triangles used in the construction of the connected Sierpi\'nski-like fractals in \S\ref{sierpinski examples}. 
In the following figures, we plot the smallest ball that intersects a certain number of children.  
The children that have non-empty intersection with the ball are colored red, while those that have empty intersection are light blue.  

 \begin{figure}[h]
		\centering
		\includegraphics[width=4in]{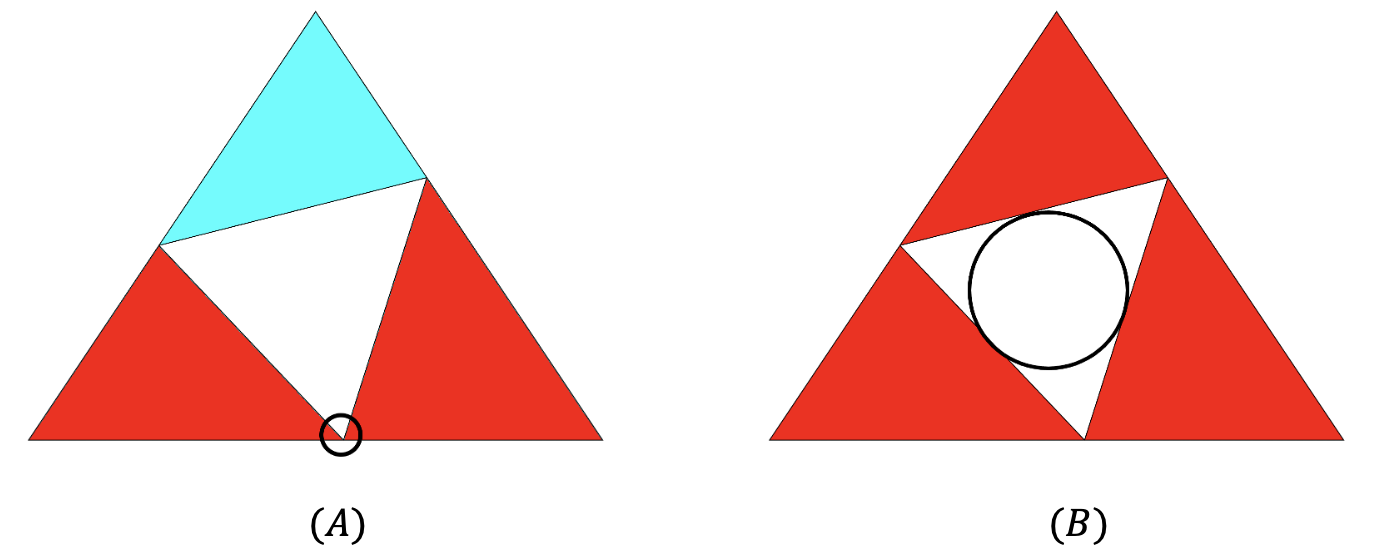}
		\caption{$(A)$ shows an illustration of $\rho_1(J_{\sigma*1}, J_{\sigma*2}, J_{\sigma*3}) = 0$ and $(B)$ shows how $\rho_2(\{J_{\sigma*1}, J_{\sigma*2}, J_{\sigma*3}\}) =$ diameter of inscribed circle}
		\label{rho1_rho2}
\end{figure}
First note that for any $\sigma \in D$, $\rho_1(\{J_{\sigma*1}, J_{\sigma*2}, J_{\sigma*3}\})=0$ since any pair of children share a vertex.  
At the intersection of the two children of $J_{\sigma}$ one can construct a ball of arbitrarily small diameter (see Figure \ref{rho1_rho2} $(A)$).
On the other hand, $\rho_2(\{J_{\sigma*1}, J_{\sigma*2}, J_{\sigma*3}\})>0$ because the diameter of any ball that intersects all three children of $J_{\sigma}$ is bounded below by the diameter of the inscribed circle of the removed center triangle.  
In other words, $\rho_2(\{J_{\sigma*1}, J_{\sigma*2}, J_{\sigma*3}\}) $ is equal to the diameter of the inscribed circle (see Figure \ref{rho1_rho2} $(B)$).

Now we may compute $\gamma_n(\mathcal{J})$ as follows. Note that for $n=1$,
\begin{align*}
		\gamma_1(\mathcal{J}) &= \inf \left \{ \frac{\rho_1(\{J_{\sigma*i}: \sigma\in R_k, i=1,2,3\})}{diam(J_\sigma)}: \text{ for some } R_k\subseteq D_k \text{ with }|R_k|= 1\right \}\\
		&= \inf \left \{ \frac{\rho_1(\{J_{\sigma*1}, J_{\sigma*2}, J_{\sigma*3}\}) }{diam(J_{\sigma})} : \text{ for } \sigma \in D \right \}=0.
	\end{align*}
On the other hand, when $n=2$, we have
\begin{eqnarray*}
	\gamma_2(\mathcal{J}) = \inf \left \{ \frac{\rho_2(\{J_{\sigma*i}: \sigma\in R_k, i=1,2,3\})}{\sum_{\sigma \in R_k}diam(J_\sigma)}: \text{ for some } R_k\subseteq D_k \text{ with }|R_k|\le 2\right \}.\\
\end{eqnarray*}
When $|R_k|=1$, this is reduced to the same case as Figure \ref{rho1_rho2} $(B)$.
When $|R_k| = 2$, we use two triangles in $R_k$, and find the smallest diameter among all balls that intersect three or more of the triangles' children. 
See Figure \ref{sierp_rho_2} for a few candidates.
\begin{figure}[h]
	\centering
	\includegraphics[scale=.08]{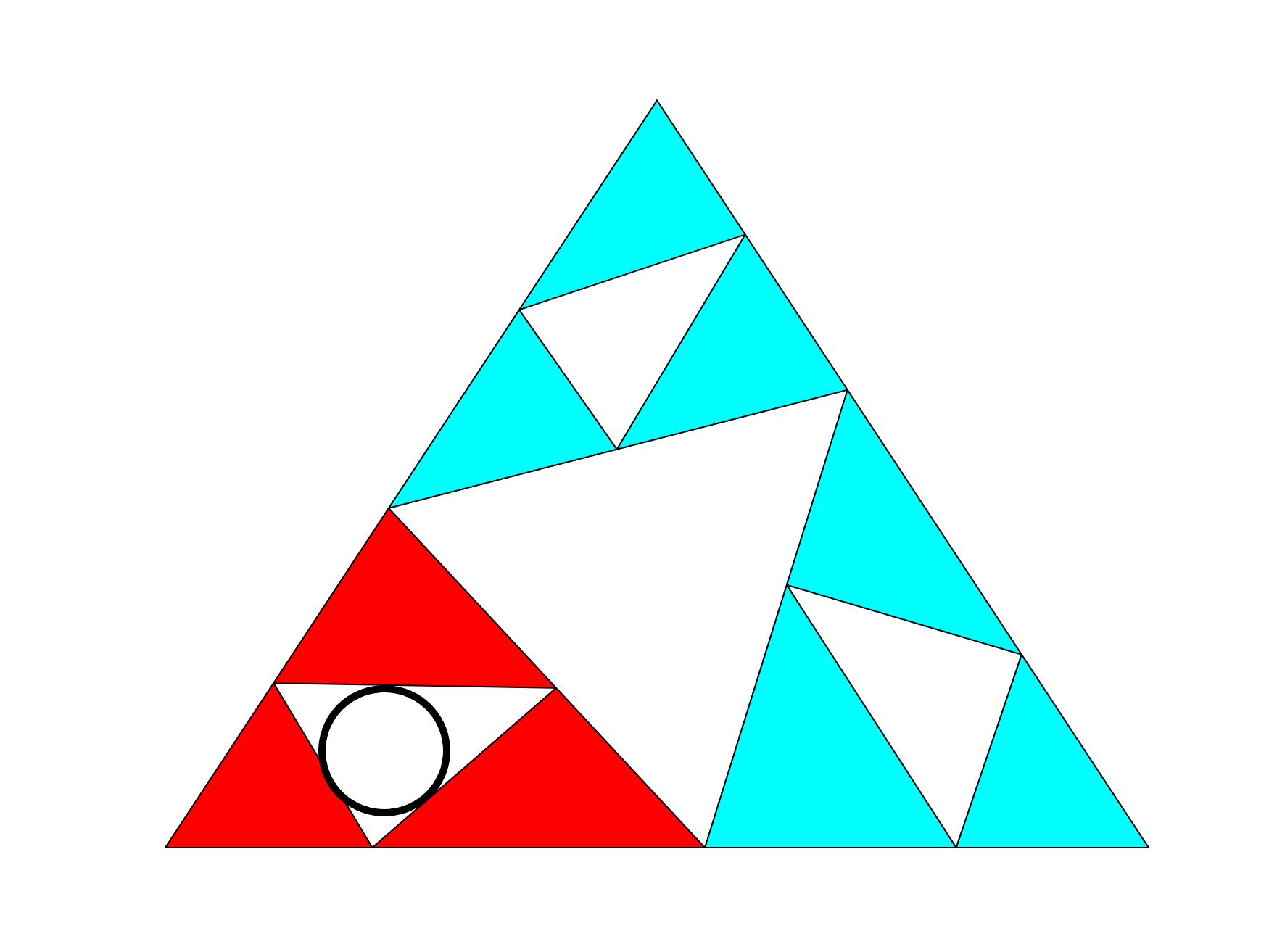}
	\includegraphics[scale=.08]{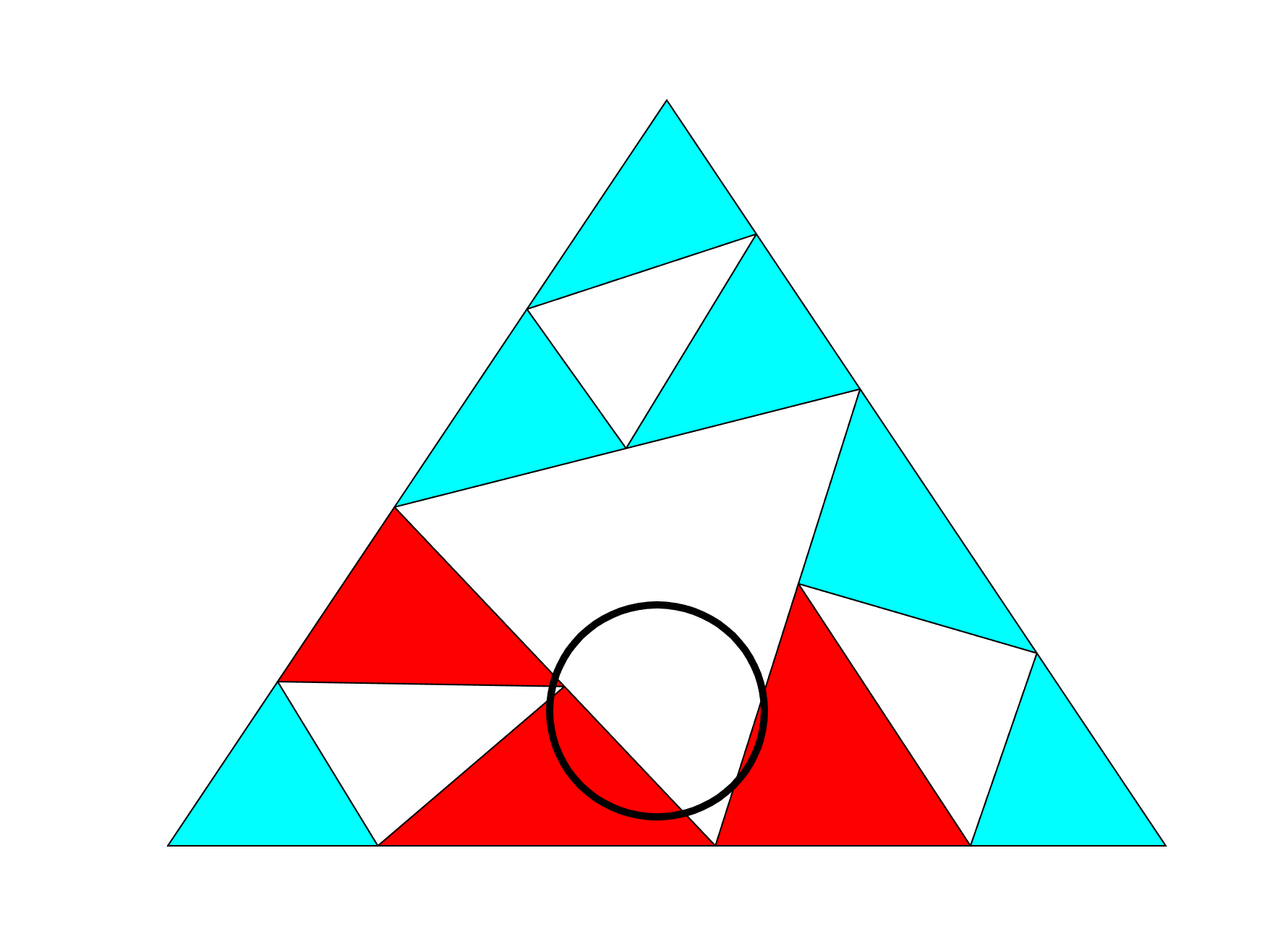}
		\includegraphics[scale=.08]{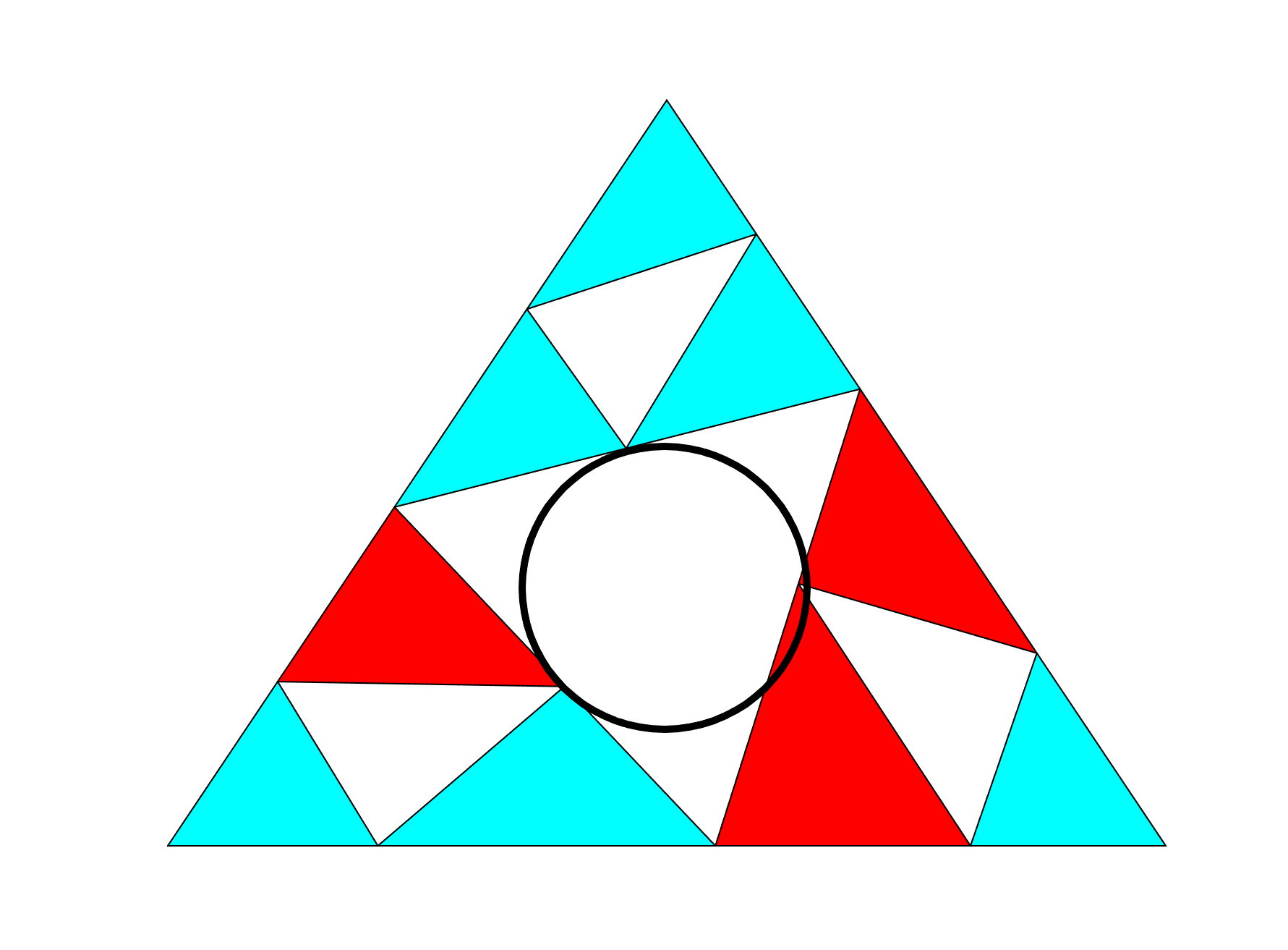}\\
	\includegraphics[scale=.08]{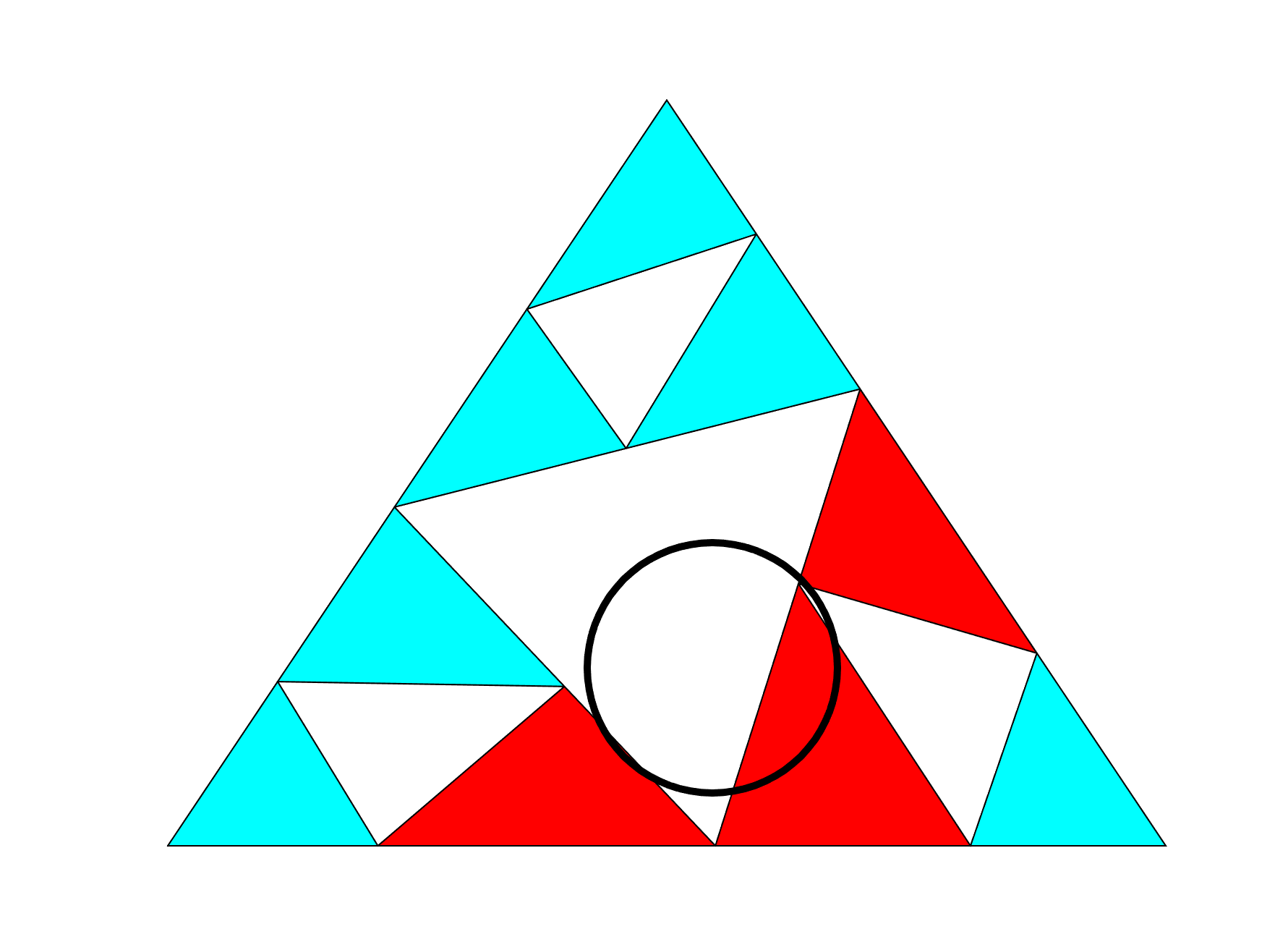}
		\includegraphics[scale=.08]{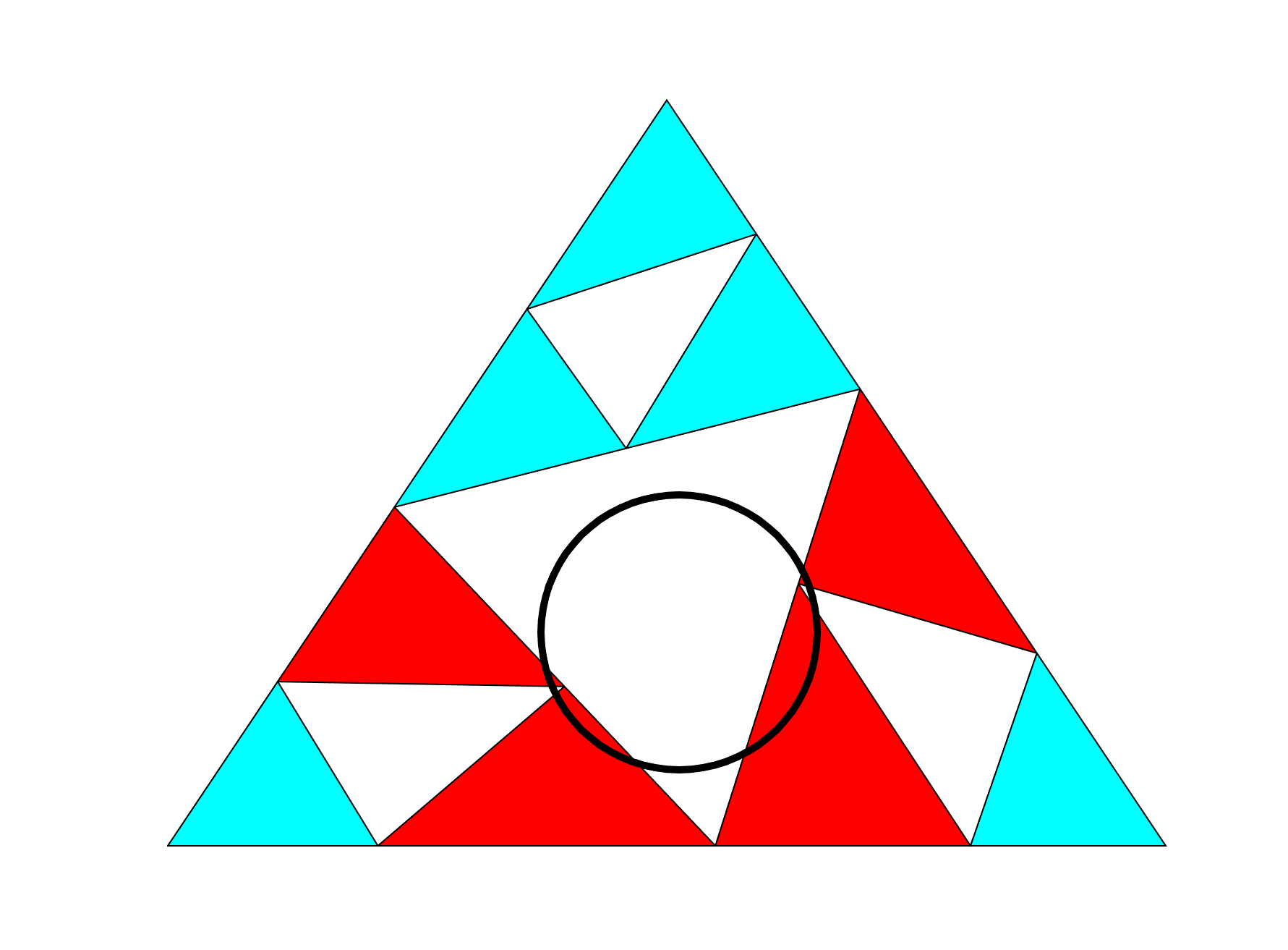}
		\includegraphics[scale=.25]{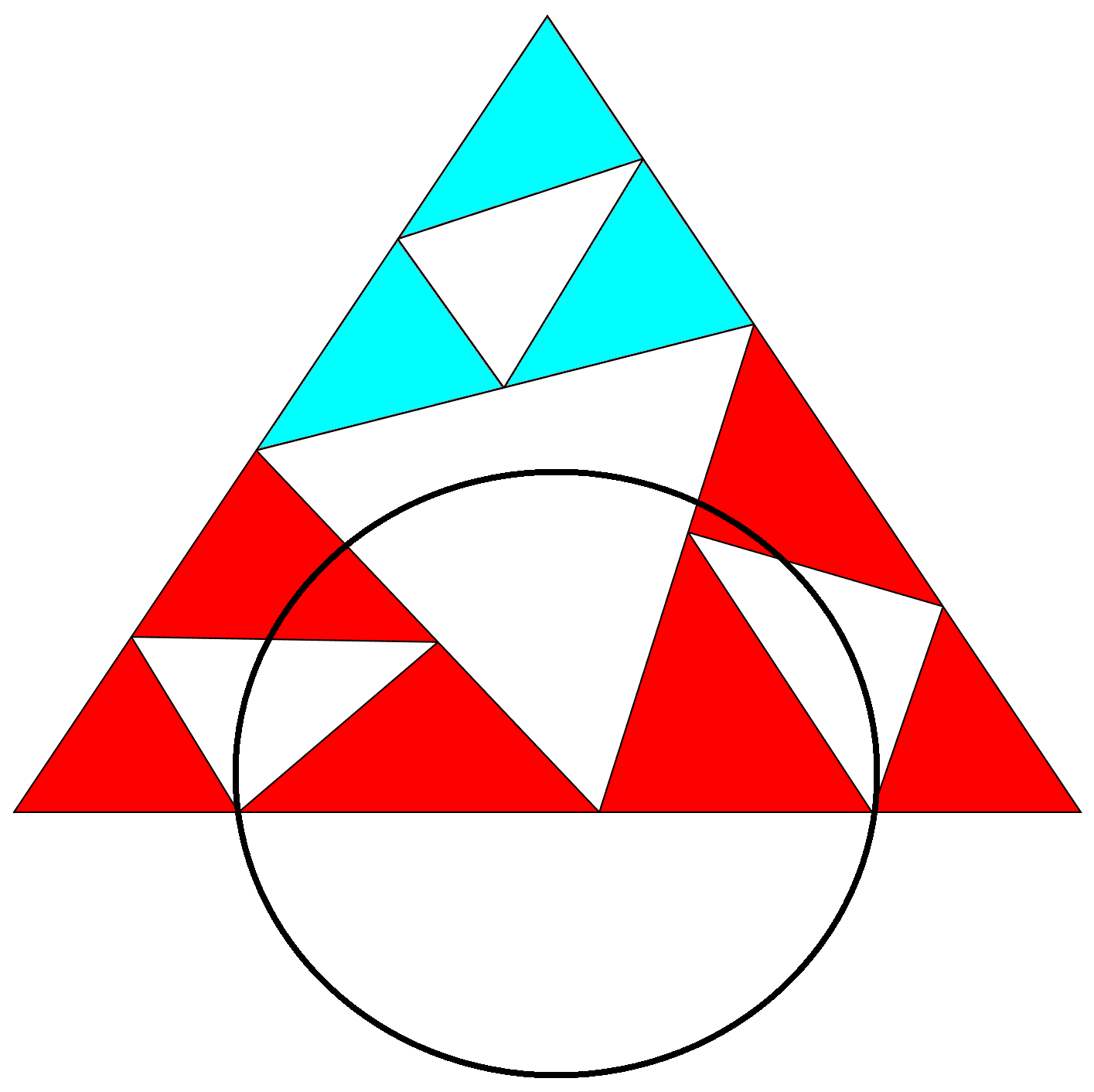}	
	\caption{Various options for smallest radius ball}
	\label{sierp_rho_2}
\end{figure}
For each $R_k \subseteq D_k$ with $|R_k| \leq 2$,  $\rho_2(\{J_{\sigma*i}: \sigma\in R_k, i=1,2,3\}) >0$. 
Note that if there exists some $c_2>0$ such that 
\[
\rho_2(\{J_{\sigma*i}: \sigma\in R_k, i=1,2,3\})\ge c_2 \max_{\sigma \in R_k} diam(J_\sigma),
\]
for all $R_k \subseteq D_k$ with $|R_k| \leq 2$ and $k = 0, 1, \dots$, then $\gamma_2(\mathcal{J})\ge c_2/2 > 0$.
\end{example}

In general, for any collection $\mathcal J := \{J_{\sigma} : \sigma \in D\}$ in a metric space $(X, d)$ and $N \in \mathbb{N}$, if there exists a constant $c_N>0$ such that
\[\rho_N(\{J_{\sigma*i}: \sigma\in R_k, i=1,2,\cdots, m\})\ge c_N \max_{\sigma \in R_k}diam(J_\sigma),\]
for any $R_k \subseteq D_k$ with $|R_k| \leq N$ and $k = 0, 1, \dots$, then $\gamma_N(\mathcal{J})\ge \frac{c_N}{N}>0$.

\begin{theorem}
	\label{thm: sufficient_uniform_covering}
	Let $\mathcal J := \{J_{\sigma} : \sigma \in D\}$ be a collection of compact subsets of $(X, d)$ satisfying MSC(3) and  
 \begin{equation}\label{eqn: max_diam}
\lim_{k\rightarrow \infty}\max \left \{diam(J_{\sigma}): \sigma\in D_k \right \}=0,
 \end{equation}
 and let $F$ be the limit set of $\mathcal{J}$ as given in (\ref{eqn: Limit_set}).  If there exists an $N$ such that $\gamma_N(\mathcal{J})>0$, then $F$ satisfies the uniform covering condition (\ref{condtion_lower}).
\end{theorem}	
\begin{proof} Let $\gamma=\gamma_N(\mathcal{J})>0$.
For any closed ball $B$ in $X$ with $B\cap F\neq \emptyset$,
let $g(k)$ be the number of elements $\sigma$ in $D_k$  such that $B\cap F\cap J_\sigma \neq \emptyset$. Then $g:\mathbb{N}\cup\{0\}\rightarrow \mathbb{N}$ is monotone increasing with $g(0)=1$.

 Case 1: Suppose $g(k)\le N$ for all $k=0,1,2,\cdots$.  For each $k$, 
 let
 \[R_k=\{\sigma\in D_k: B\cap F\cap J_\sigma \neq \emptyset\}\subseteq D_k.\]
 Then, $|R_k|=g(k)\le N$  and $B\cap F\subseteq \bigcup_{\sigma \in R_k}J_{\sigma}$. As a result,
 \[0\le \sum_{\sigma \in R_k}diam(J_{\sigma}) \le N\cdot \max\{diam(J_{\sigma}): \sigma\in D_k\}.\]
By (\ref{eqn: max_diam}) and the squeeze theorem,
\[\lim_{k\rightarrow \infty} \sum_{\sigma \in R_k}diam(J_{\sigma})=0.\]
Thus, since $\gamma>0$, when $k$ is large enough,
\[diam(B)>\gamma \cdot \sum_{\sigma\in R_k}diam(J_{\sigma}).\]
 Hence, equation  (\ref{condtion_lower}) holds for $B$.

Case 2: There exists $k^* \geq 0$ such that $g(k^*)\leq N$ but $g(k^*+1)>N$.
	
Since $g(k^*) \leq N$, there are $g(k^*)$ many elements $\sigma \in D_{k^*}$ such that $B \cap F \cap J_{\sigma} \neq \emptyset.$  That is, there exists $R_{k^*} \subseteq D_{k^*}$ with $|R_{k^*}| \leq N$ such that $B \cap F \subseteq \bigcup_{\sigma \in R_{k^*}} J_{\sigma}$.  On the other hand, since $g(k^*+1)>N$, $B \cap F$ intersects at least $N+1$ elements of $D_{k^*+1}$.  Since $B\cap F \subseteq  \bigcup_{\sigma \in R_{k^*}} J_{\sigma}$, all of these $N+1$ elements must be children of $\{J_{\sigma} : \sigma \in R_{k^*}\}$.  Then, by the definition of $\rho_N$ in (\ref{rho_n}),  \begin{equation}
	diam(B) \geq \rho_N( \{J_{\sigma*i}: \sigma \in R_{k^*}, i=1,2,\dots m\}) \geq \gamma \cdot \sum_{\sigma \in R_{k^*}} diam(J_{\sigma}).	\end{equation}
	 As a result, $F$ satisfies the uniform covering condition (\ref{condtion_lower}). 
  \end{proof}

  Now we explicitly apply Theorem \ref{thm: sufficient_uniform_covering} to the Menger sponge-like fractal of Example~\ref{3d_exact}.
  Let $\mathcal{X}$, $D$, $M$, and $\mathcal{F}$ be as in section \S\ref{menger section}.
  To show that the $\mathcal{F}$-limit set generated by certain sequences satisfies the uniform covering conditions (\ref{condtion_lower}), we will first derive the following lower bound~(\ref{lower bound}) of $\gamma_8(\mathcal{J})$ for a collection $\mathcal{J} = \{J_\sigma : \sigma \in D\}$ generated by an {\it arbitrary} sequence $\{\k_\ell\}_{\ell = 0}^\infty$ in $M$. 
	
Let $\mathcal{H}$ be a subset of $\{J_\sigma : \sigma \in D_k\}$ for some $k \in \{0, 1, \dots\}$. 
We now make the following observation: \textit{For any ball $B$ that intersects at least 9 elements of $\mathcal{H}$, its diameter  $diam(B)$ must be greater than or equal to the smallest edge length of the elements in $\mathcal{H}$.}  
Indeed, by considering the projections to the three coordinate axes, one can see that at least one coordinate contains three non-identical projected images of these 9 elements.  
As a result, the ball $B$ intersected with these 9 elements will have a diameter at least the length of the smallest side of the three projected images.  
This justifies our observation. 

For each $\k_\ell = (k_\ell^{(1)}, k_\ell^{(2)},k_\ell^{(3)},k_\ell^{(4)},k_\ell^{(5)},k_\ell^{(6)})  \in M$, define 
\begin{equation*}
	m_{\ell} = \min \{ k_\ell^{(1)}, k_\ell^{(2)}-k_\ell^{(1)},1-k_\ell^{(2)}, k_\ell^{(3)}, k_\ell^{(4)}-k_\ell^{(3)},1-k_\ell^{(4)}, k_\ell^{(5)}, k_\ell^{(6)}-k_\ell^{(5)},1-k_\ell^{(6)} \}
	\end{equation*} and 
	\begin{equation*}
	M_{\ell} = \max \{ k_\ell^{(1)}, k_\ell^{(2)}-k_\ell^{(1)},1-k_\ell^{(2)}, k_\ell^{(3)}, k_\ell^{(4)}-k_\ell^{(3)},1-k_\ell^{(4)}, k_\ell^{(5)}, k_\ell^{(6)}-k_\ell^{(5)},1-k_\ell^{(6)} \}.
\end{equation*}
Then, by the definition of $M$, we have $0\le m_\ell\le M_\ell\le 1$ for each $\ell$.

For any $\sigma  \in D$, direct calculation shows that \begin{equation}
m_{\ell(\sigma)} \leq \dfrac{diam(J_{\sigma*i})}{diam(J_{\sigma})} \leq M_{\ell(\sigma)}
\end{equation}
where $\ell(\sigma)$ is given in (\ref{ell_sigma}). Thus, for any $\sigma = (i_1, i_2, \dots, i_{k}) \in D_{k}$, we have \begin{equation}
	m_{\ell((i_1))} m_{\ell((i_1,i_2))} \cdots m_{\ell((i_1,\dots,i_k))} \leq \dfrac{diam(J_{\sigma})}{diam(J_{\emptyset})} \leq 	M_{\ell((i_1))} M_{\ell((i_1,i_2))} \cdots M_{\ell((i_1,\dots,i_k))}.
	\end{equation}
	
Let $R_{k} \subseteq D_{k}$ for some $k$.  Suppose $|R_{k}|\leq8$.  Then for any $\sigma \in R_{k}$, by the observation \begin{eqnarray*}
	&& \dfrac{\rho_8(\{J_{\sigma*i}: \sigma \in R_{k} , i=1,2, \dots 20\})}{\sum_{\sigma \in R_k} diam(J_{\sigma})}  \\
	&\geq&
	\dfrac{ \text{ smallest diameter of $J_{\sigma*i}$}}{8 \cdot \max \{ diam(J_{\sigma}): \sigma \in R_{k}\}}\\
	&\geq& \dfrac{1}{8} \min_{(i_1, i_2, \dots , i_k) \in R_k, i_{k+1}=1, \dots , 20} \left \{ \dfrac{m_{\ell((i_1))} m_{\ell((i_1,i_2))} \cdots m_{\ell((i_1,\dots,i_{k+1}))}  diam(J_{\emptyset})}{M_{\ell((i_1))} M_{\ell((i_1,i_2))} \cdots M_{\ell((i_1,\dots,i_k))} diam(J_{\emptyset})}  \right\}\\
	&\geq& \dfrac{1}{8} \left (\prod_{i=1}^{\infty} \dfrac{ m_{i}}{M_{i}} \right) \liminf_{i \to \infty} m_{i},\\
\end{eqnarray*}  
where the last inequality follows from $0 \leq m_i \leq M_i$ for each $i$.	As a result, we have
\begin{equation}\label{lower bound}
\gamma_8(\mathcal{J})\ge \dfrac{1}{8} \left (\prod_{i=1}^{\infty} \dfrac{ m_{i}}{M_{i}} \right) \liminf_{i \to \infty} m_{i}.
\end{equation}
	
\begin{example}\label{3d_ucc}
We now apply it to show that the $\F$-limit set in Example \ref{3d_exact} satisfies the uniform covering condition.	
In this example, 
\begin{equation*}
	 m_{\ell} = \begin{cases}
	 a_{\ell}, & \ell \text{ even}\\
	 b_{\ell}, & \ell \text{ odd}
	 \end{cases}
	 \text{ and }
	 M_{\ell} = \begin{cases}
	 b_{\ell}, & \ell \text{ even}\\
	 a_{\ell}, & \ell \text{ odd}
	 \end{cases}
\end{equation*}
where
 \begin{equation*}
	a_{\ell} = k_{\ell}^{(3)}= \frac{1}{3} - \frac{(-1)^{\ell}}{6(\ell +1)^2} \text{ and } b_{\ell}= 1-2 k_{\ell}^{(3)}= \frac{1}{3} + \frac{(-1)^{\ell}}{3(\ell+1)^2}.
	\end{equation*}  
	
One may show that the product $\displaystyle \prod_{i=1}^{\infty} \dfrac{ m_{i}}{M_{i}} $ is convergent, whose numerical value is 0.369761$\dots$ and $\liminf_{i \to \infty} m_i =1/3$.  Thus, by (~\ref{lower bound}), $\gamma_8(\mathcal{J}) >0$.  Therefore, by Theorem \ref{thm: sufficient_uniform_covering}, the $\mathcal{F}$-limit set $F$ satisfies the uniform covering condition.
\end{example}

\end{document}